\documentclass[11pt]{article}
\usepackage{epsfig,amsmath,latexsym}
\usepackage{amsfonts}
\newtheorem{theorem}{Theorem}[section]

\newtheorem{defi}{Definition}[section]
\newtheorem{lemma}{Lemma}[section]

\def\binom#1#2{{#1}\choose{#2}}

\def\slfrac#1#2{\hbox{\kern.1em %
 \raise.5ex\hbox{\the\scriptfont0 #1}\kern-.11em %
 /\kern-.15em\lower.25ex\hbox{\the\scriptfont0 #2}}}

\newcommand{\eqn}[1]{(\ref{#1})}
\newcommand{\hsp}{\hspace*{\parindent}}

\newcommand{\eeq}{\end{equation}}
\newcommand{\beql}[1]{\begin{equation}\label{#1}}
\newcommand{\bsq}{{\vrule height .9ex width .8ex depth -.1ex }}

\newcommand{\AAA}{{\Bbb A}}
\newcommand{\CC}{{\Bbb C}}
\newcommand{\FF}{{\Bbb F}}

\newcommand{\QQ}{{\Bbb Q}}
\newcommand{\RR}{{\Bbb R}}
\newcommand{\ZZ}{{\Bbb Z}}

\newcommand{\sA}{{\cal A}}

\newcommand{\sH}{{\cal H}}
\newcommand{\sL}{{\cal L}}
\newcommand{\sS}{{\cal S}}
\newcommand{\sW}{{\cal W}}
\newcommand{\fq}{{\mathfrak q}}

\makeatletter
\def\@sect#1#2#3#4#5#6[#7]#8{\ifnum #2>\c@secnumdepth
     \def\@svsec{}\else
     \refstepcounter{#1}\edef\@svsec{\csname the#1\endcsname.\hskip .75em }\fi
     \@tempskipa #5\relax
      \ifdim \@tempskipa>\z@
        \begingroup #6\relax
          \@hangfrom{\hskip #3\relax\@svsec}{\interlinepenalty \@M #8\par}%
        \endgroup
       \csname #1mark\endcsname{#7}\addcontentsline
         {toc}{#1}{\ifnum #2>\c@secnumdepth \else
                      \protect\numberline{\csname the#1\endcsname}\fi
                    #7}\else
        \def\@svsechd{#6\hskip #3\@svsec #8\csname #1mark\endcsname
                      {#7}\addcontentsline
                           {toc}{#1}{\ifnum #2>\c@secnumdepth \else
                             \protect\numberline{\csname the#1\endcsname}\fi
                       #7}}\fi
     \@xsect{#5}}
\def\@begintheorem#1#2{\it \trivlist \item[\hskip \labelsep{\bf #1\ #2.}]}

\def\plain{plain}\ifx\fmtname\plain\csname fi\endcsname
     
     \input docstrip
     \preamble

     Do not distribute the stripped version of this file.
     The checksum in the header refers to the documented version.

     \endpreamble
     \generateFile{here.sty}{t}{\from{here.doc}{}}
     \endinput
\fi
\ifcat a\noexpand @\let\next\relax\else\def\next{%
    \documentstyle[here,doc]{article}\MakePercentIgnore}\fi\next
\ifx\@Hxfloat\@Hundef\else\expandafter\endinput\fi
\let\@Hxfloat\@xfloat
\def\@xfloat#1[{\@ifnextchar{H}{\@HHfloat{#1}[}{\@Hxfloat{#1}[}}
\def\@HHfloat#1[H]{%
\expandafter\let\csname end#1\endcsname\end@Hfloat
\vskip\intextsep\vbox\bgroup\def\@captype{#1}\parindent\z@
\ignorespaces}
\def\end@Hfloat{\egroup\vskip \intextsep}

\makeatother

\catcode`\@=11
\renewcommand{\section}{
        \setcounter{equation}{0}
        \@startsection {section}{1}{\z@}{-3.5ex plus -1ex minus
        -.2ex}{2.3ex plus .2ex}{\large\bf}%
        }
\catcode`@=12

\begin{document}
\begin{center}
{\Large {\bf  Li  Coefficients for Automorphic $L$-Functions}}\\

\vspace{1.5\baselineskip}
{\em Jeffrey C. Lagarias} \\
\vspace*{.2\baselineskip}
University of Michigan  \\
Ann Arbor, MI 48109-1043\\
email: {\tt lagarias@umich.edu} \\
\vspace*{2\baselineskip}
(April 11, 2005) \\
\vspace{1.5\baselineskip}
{\bf Abstract}
\end{center}
\noindent
Xian-Jin Li gave a criterion for the Riemann
hypothesis in terms of the positivity of
the set of coefficients
$\lambda_n = \sum_{\rho} 1 - \left( 1 - \frac{1}{\rho}\right)^n$,
$(n= 1, 2, ...)$,
in which $\rho$ runs over the nontrivial zeros of the
Riemann zeta function. We define  similar 
coefficients $\lambda_n(\pi)$ associated to principal automorphic
$L$-functions $L(s, \pi)$ over $GL(N)$.
 We relate these cofficients
to values of Weil's quadratic functional 
associated to the representation $\pi$ on a suitable
set of test functions. The positivity of the
real parts of these coefficients is a necessary and sufficient 
condition for the Riemann hypothesis for $L(s, \pi)$ to hold.

We derive an unconditional
asymptotic  formula for the coefficients $\lambda_n(\pi)$,
in terms of the zeros of  $L(s, \pi)$.
Assuming  the Riemann hypothesis for $L(s, \pi)$, we deduce 
that \\
$\lambda_n(\pi)  = \frac{N}{2} n \log n + C_1(\pi) n 
+ O (\sqrt{n}\log{n}),$
where $C_1(\pi)$ is a real-valued constant
and the implied constant in the remainder term depends on $\pi$.
We also show that  there exists a
entire function $F_{\pi}(z)$ of
exponential type  that 
interpolates the generalized Li coefficients at integer values.
Assuming the Riemann hypothesis there is an (essentially) 
unique interpolation function having  exponential
type at most $\pi$, and this function restricted to the
real axis has a (tempered) distributional Fourier transform
whose support is a countable
set in $[-\pi, \pi]$ having $0$ as its only limit point.

 
%
%
%
\setlength{\baselineskip}{1.0\baselineskip}

\section{Introduction}
\hsp
In 1997 Xian-Jin Li \cite{Li97} derived a necessary and
sufficient condition for the Riemann hypothesis in
terms of the positivity of the set of coefficients
\beql{103}
\lambda_n := \sum_{\rho} {}^{'} [1 - \left( 1 - \frac{1}{\rho}\right)^n], 
\eeq 
in which the sum runs over the nontrivial
zeros of the Riemann zeta function, counted with
multiplicity, and ${}^{'}$ indicates that the 
(conditionally convergent) sum is to be
interpreted as
$\lim_{T \to \infty} \sum_{ \{\rho:~ |\rho| \le T \} };$
we term this  $*$-convergence, 
The  expression \eqn{103} $*$-converges for 
positive and negative integer $n$, and
so defines these coefficients for all integers.
One has $\lambda_0 = 0$ and 
\begin{equation}~\label{103c}
\lambda_{-n} = \lambda_{n},~~  \mbox{for~ all}~~ n \ge 1.
\end{equation}
Li's original definition
\beql{101}
\tilde{\lambda}_n := 
\frac{1}{(n-1)!}\frac{d^n}{ds^n} [s^{n-1}  \log \xi(s)] |_{s=1}
~~~ n \ge 1,
\eeq
functorially
corresponds to  $\tilde{\lambda}_n = \lambda_{-n}$ for $n > 0$,
and the identity
$\tilde{\lambda}_{n} =\lambda_{n}$ then holds using \eqn{103c}. 
Li's paper gave more generally a criterion for 
the Riemann hypothesis to hold  for the (completed)
Dedekind zeta functions of any algebraic number field $K$.

These coefficients are expressible in
terms of power-series coefficients of functions constructed
from the Riemann $\xi$-function,
$\xi(s) = \frac{1}{2} s(s-1)\pi^{-s/2} \Gamma(\frac{s}{2}) \zeta(s)$.
The coefficients for $n \ge 1$ occur in 
\beql{102}
\frac{d}{dz} \log \xi(\frac{-z}{1-z}) =
\frac{-1}{(1-z)^2}\frac{\xi'}{\xi}(\frac{-z}{1-z}) = 
\sum_{n=0}^{\infty} \lambda_{n+1} z^n,
\eeq 
and those for $n \le -1$ in occur in   an analogous formula
\beql{102a}
\frac{d}{dz} \log \xi(\frac{1}{1-z}) =
\frac{1}{(1-z)^2}\frac{\xi'}{\xi}(\frac{1}{1-z}) = 
\sum_{n=0}^{\infty} \lambda_{-n-1} z^n.
\eeq
 Here $\xi'(\cdot)$ denotes the derivative with respect to $s$.
Our definition \eqn{103} of the Li coefficients 
corresponds to an expansion around the point $s=0$,  visible in \eqn{102},
rather than around $s=1$ as in \eqn{101} and \eqn{102a}.

Closely related coefficients
appear in earlier work of J. B. Keiper \cite[Sect. 4]{Ke92} in 1992, 
in an investigation of methods to compute Stieltjes constants to
high precision; we describe Keiper's computations 
futher below.  Keiper's coefficients
equal $\frac{1}{n}\lambda_{-n}$, in terms of  \eqn{103}.
Keiper noted that the Riemann hypothesis implies
the nonnegativity of the $\lambda_{-n}$, but not the converse. 
We adopt the term 
{\em Li coefficients} because Li's work shows the nonnegativity of 
$\lambda_n$ is equivalent to the Riemann hypothesis and 
because he generalized them to various other zeta functions.

In \cite{BL99} E. Bombieri and the author
made three observations about these coefficients.  The first 
observation was  
that a Li criterion can be formulated for 
very general sets of complex numbers $\rho$,
as follows.  Consider any
multiset $Z$ of complex numbers $\rho$ satisfying
\beql{102b}
\sum_{\rho \in Z} \frac{ \Re(\rho)}{ (1 + |\rho|)^2} < \infty.
\eeq
If the multiset $Z$ omits the value $\rho=1$ then 
the sums 
\beql{103a}
\Re(\lambda_n(Z)) : = \sum_{\rho \in Z} 
\Re\left( 1 - ( 1- \frac{1}{\rho})^n \right)
\eeq
converge absolutely  for all 
nonpositive integers $n \le 0$. The
positivity condition $\Re(\lambda_n(Z))\ge 0 $ for $n \le 0$
then implies that all $\rho$ lie in the half-plane 
$\Re(s) \le \frac{1}{2}.$ 
If the multiset $Z$ omits the value $\rho=0$, then 
the sum \eqn{103a} then converges absolutely for $n \ge 0$,
and the positivity condition $\Re(\lambda_n(Z))\ge 0 $ for $n \ge 0$
implies 
that all $\rho$ lie in the half-plane $\Re(s) \ge \frac{1}{2}.$
Combining these criteria, for 
multisets $Z$ that omit both $0$ and $1$ the positivity
condition  $\Re(\lambda_n(Z)) \ge 0 $ for all integers $n$
implies that all  $\Re(\rho) = \frac{1}{2}$.
If the multiset $Z$ is also
invariant under the symmetry 
$\rho \mapsto 1 - \bar{\rho}$,
so that  $\Re(\lambda_n(Z)) = \Re(\lambda_{-n}(Z))$,
it suffices to check 
this positivity condition $\Re(\lambda_n(Z)) \ge 0$ for  $n > 0$
to conclude that all  $\Re(\rho) = \frac{1}{2}$. 
Finally, if  $Z$  omits the values $0$ and $1$ and
the sum $\sum_{\rho \in Z} \frac{1}{\rho}$ 
is $*$-convergent, then the coefficients $\lambda_n(Z)$
are well-defined for all integers $n$ by the following 
$*$-convergent sum:
\begin{equation}~\label{103e}
\lambda_n(Z) := \sum_{\rho \in Z}{}^{'} 
\left( 1 - ( 1- \frac{1}{\rho})^n \right).
\end{equation}
The second observation in \cite[Theorem 2]{BL99} was that  the
the ``explicit formula'' of prime number theory may be used to obtain 
an arithmetic expression for Li's coefficients $\lambda_n$
in \eqn{103} having the form 
$$
\lambda_n = S_{\infty}(n) - S_{f}(n) + 1,
$$
in which  $S_{\infty}(n)$ and $S_{f}(n)$ correspond 
to the contributions of the archimedean place and the finite
places, respectively, and the last term is a contribution
from the pole at $s=0$ of $\frac{\xi(s)}{s(s-1)}$. 
The third observation
was that  each positivity condition  $\lambda_n \ge 0$  
encodes ``Weil positivity''
of Weil's quadratic functional 
for a particular test function $g_n(x)$.

In more recent work K. Maslanka~\cite{Ma03} 
computed $\lambda_n$ 
for $1 \le n \le 3300$ and 
empirically studied
the growth behavior of the Li coefficients. 
He observed that in this range they
exhibited  a smoothly growing 
dominant asymptotic term 
with superposed  small oscillations exhibiting some internal structure.
The  dominant term comes  
from the archimedean prime contribution $S_{\infty}(n)$ 
in the arithmetic formula above, and the small oscillations
come  from the term $S_{f}(n)$ represented the finite places,
and $|S_{f}(n)| < 20$ over this range.
Coffey \cite{Co03},\cite{Co03b} studied the ``arithmetic formula'' for
the Li coefficients and lower  bounded
the archimedean prime contribution.

The object of this  paper is to generalize 
the Li coefficients to automorphic $L$-functions,
and to determine their asymptotic behavior as $n \to \infty$,
with or without the assumption of the Riemann hypothesis.
The automorphic $L$-functions
we treat are  principal $L$-functions over $GL(N)$
for $\QQ$, as given in \cite{Cog03}, \cite{IK04},
\cite{Ja79}, and \cite{RS96}.
These are (completed) Langlands $L$-functions  attached
to irreducible cuspidal unitary automorphic representations 
appearing in the right action of $GL(n, \AAA_{\QQ})$ on
$L^2(GL(n, \QQ) \backslash GL\left(n, \AAA_{\QQ})\right)$.
The associated 
generalized Li coefficients will be denoted $\lambda_n(\pi).$

In \S2 we associate to each
irreducible cuspidal  unitary automorphic representation $\pi$  
an analogue $\xi(s, \pi)$ of 
the Riemann $\xi$-function.
Our definition of  {\em generalized  Li coefficients} $\lambda_n(\pi)$ 
sets them equal to $\lambda_n(Z)$ for the multiset
$Z= Z(\pi)$ of zeros of  $\xi(s, \pi)$,
and we show that these coefficients $\lambda_n(\pi)$ are
well-defined as $*$-convergent series.
The multiset $Z(\pi)$ is invariant under
the transformation $\rho \mapsto 1 - \bar{\rho}$, which
implies the symmetry 
\beql{107a}
\lambda_{-n}(\pi) = \overline{ \lambda_n(\pi) }.
\eeq
It therefore suffices to study $\lambda_n(\pi)$ for $n \ge 0$.
Furthermore the zero sets of the $L$-functions
 of a representation $\pi$ and its contragredient representation
$\pi^{\vee}$ are
related by complex conjugation $Z(\pi^{\vee}) = \overline{Z(\pi)}$,
which implies a second symmetry
\beql{107b}
\lambda_{n}(\pi^{\vee}) = \overline{\lambda_{n} (\pi)}.
\eeq
The results of \cite{BL99}  apply to give a
Riemann hypothesis criterion for $L(s, \pi)$ in the
form: The Riemann hypothesis holds for $L(s, \pi)$ 
if and only if the real parts  $\Re(\lambda_n(\pi))$ are nonnegative for all 
 $n \ge 0.$

In \S3 we  give a ``Weil positivity'' interpretation of
the generalized Li coefficients. We express the  Weil scalar
product associated to 
the representation $\pi$ for Li's test functions
in terms of $\lambda_n(\pi)$. Here  the imaginary 
parts of the $\lambda_n(\pi)$ appear in the scalar products.
The test functions are the same for all representations $\pi$.

In \S4 we give an arithmetical
interpretation of the coefficients in
terms of  the logarithmic
derivative of $\xi(s, \pi)$ expanded about at the point $s=1$.
We express it as 
\beql{105b}
 \lambda_{n}(\pi)=
  S_{\infty}(n, \pi^{\vee}) - S_{f}(n, \pi^{\vee}) + 
\delta(\pi^{\vee}),
\eeq
in which the  two terms  $S_{\infty}(n, \pi^{\vee})$ and 
$S_{f}(n, \pi^{\vee})$
reflect contributions coming from  the Euler product factorization
of $\xi(s, \pi^{\vee})$ 
into archimedean places and finite places, respectively.
The final term  $\delta(\pi)=\delta(\pi^{\vee})= 1$ 
for the trivial representation $\pi= \pi_{triv}$
over $GL(1)$ and $\delta(\pi)=0$ otherwise.
The contragredient representation appears in  \eqn{105b}
via the functional equation relating $\xi(s, \pi)$
and $\xi(1-s, \pi^{\vee})$, because our definition \eqn{103}
of the Li coefficients really corresponds to values at $s=0$.

In \S5 we obtain an unconditional asymptotic formula  for 
the archimedean contribution $S_{\infty}(n, \pi)$.
This quantity is real-valued, and we show
there is a real-valued constant $C_1(\pi)$
such that for all $n \ge 1$, 
\beql{105c}
S_{\infty}(n, \pi) = \frac{N}{2} n \log n + C_1(\pi)~ n + 
O \left( 1\right),
\eeq
and  the implied constant 
in the $O(1)$ term depends on $\pi$.  Here
\beql{105d}
C_1(\pi) = \frac{N}{2}( \gamma - 1 - \log(2\pi)) + \frac{1}{2} \log  Q(\pi),
\eeq
in which $\gamma$ is Euler's constant
and $Q(\pi)$ is the conductor of $\pi$. In particular
$C_1(\pi)$
does not depend on  
the archimedean parameters
$\{ \kappa_j(\pi): 1 \le j \le N\}$ of the representation
$\pi$, and $C_1(\pi)= C_1(\pi^{\vee})$. 

In \S6 we obtain an unconditional estimate for 
the finite place contribution $S_{f}(n, \pi)$
in terms of the zeros to a suitable height, by a contour
integral estimate.  
Define the {\em incomplete  Li coefficient} at height $T$ by 
\beql{106}
\lambda_n(T, \pi) = \sum_{{{\rho \in Z(\pi)}\atop{|\Im(\rho)| < T}}}
1 - \left( 1 - \frac{1}{\rho} \right)^n,
\eeq
We show that
\beql{107}
S_f(n , \pi) = \lambda_n(\sqrt{n}, \pi^{\vee}) + 
O \left( \sqrt{n} \log n \right),
\eeq
where $\pi^{\vee}$ is the contragredient representation, and 
the implied constant in the $O$-notation
depends on $\pi$.
If the Riemann hypothesis 
holds for $L(s, \pi)$, then it holds   for $L(s, \pi^{\vee})$
by the functional equation, and we obtain
$$
  \lambda_n(\sqrt{n}, \pi^{\vee})= O\left(\sqrt{n}\log n \right),
$$
so that 
 \beql{N107c}
S_{f}(n , \pi) = O\left( \sqrt{n} \log n \right). 
\eeq
Furthermore if the Riemann hypothesis 
holds up to height $T$, then a  bound of shape \eqn{N107c} holds for all 
$n \le \frac{T^2}{4(\log T)^2}$, with the
implied O-constant depending on $\pi$.

These results are summarized in the following theorem.

\begin{theorem}~\label{th11}
Let $\pi$ be an irreducible cuspidal 
unitary automorphic representation
for $GL(N)$ over $\QQ$. 
For $n \ge 1$ there holds
\beql{108}
\lambda_n(\pi) = 
\frac{N}{2} n \log n +  C_1(\pi) n -
\lambda_n(\sqrt{ n}, \pi) + O (\sqrt{n}\log n),
\eeq
in which $C_1(\pi)$ is real-valued and 
the implied  constant in the $O$-notation depends on $\pi$.
If the Riemann hypothesis holds  for $L(s, \pi)$ then
the incomplete Li coefficient 
$\lambda_n(\sqrt{n}, \pi) = O (\sqrt{n}\log n)$, so that 
for $n \ge 1$, 
\beql{108b}
\lambda_n(\pi) = 
\frac{N}{2} n \log n +  C_1(\pi) n +  O (\sqrt{n}\log n),
\eeq
where the implied constant in the $O$-notation depends on $\pi$. 
\end{theorem}

 Theorem 1.1 follows on combining  Lemma~\ref{le42}, Theorem~\ref{th51}
and Theorem~\ref{th61}, together with 
the relation $C_1(\pi)=C_1(\pi^{\vee})$.
If the Riemann hypothesis does not hold for $L(s, \pi)$
then the incomplete Li coefficient
term  $\lambda_n(\sqrt{n}, \pi)$ will sometimes be very large,
of size exponential in $n$.
This fact was already observed for the Riemann zeta
function in \cite[Theorem 1(c)]{BL99}.

In \S7 we construct  for each $\pi$   
an entire function $F_{\pi}(z)$  of order one
and exponential type 
that interpolates the Li coefficients at integer values,
i.e. $F_{\pi}(n) = \lambda_n(\pi)$ for all $n \in \ZZ$.
Assuming the Riemann hypothesis for $L(s, \pi)$, 
this function $F_{\pi}(z)$ can be chosen to have
exponential type at most $\pi$. 
It is then (almost) uniquely
characterized by the interpolation property;
the remaining ambiguity concerns zeros at the
central critical value $s= \frac{1}{2}$,
which require special treatment, see 
 Theorem~\ref{th71}.
The polynomial growth of the Li coefficients
under the RH allows one to deduce that the
function $F_{\pi}(z)$
has a Fourier transform which is well-defined
as a tempered distribution, whose support 
is a countable closed set in $[-\pi, \pi]$
having $0$ as its only limit point. 
This interpolation function $F_{\pi}(z)$ 
appears to be a new object associated
to the zeta zeros.  We do not know any relation
of it to various other functions constructed
from the zeta zeros such as  
Cramer's $V$-function
(\cite{Cr19}, \cite{Gu49}, \cite{DS95}, \cite{JL01}, \cite{Il02})
or functions studied by 
Voros \cite{Vo92}, \cite{Vo03}.

The results of this paper explain some of the 
empirical observations of  
K. Maslanka \cite{Ma03}.
The asymptotic formula  \eqn{N107c} applied to $\pi_{triv}$,
explains the observed behavior of the small Li coefficients.
Since the
non-trivial zeros of $\zeta(s)$ are known
to  lie on  the critical line
up to height $T \approx 10^9$ we may expect the 
first $10^{16}$ Li coefficients will also  exhibit similar 
asymptotic behavior, i.e. the term $|S_{f}(n)|$
will remain small over this range. 
Maslanka's computations of $S_{f}(n)$
allow the possibility the  term  $|S_{f}(n)|$ is 
of smaller order of growth than $O(\sqrt{n}\log n)$;
if so, this remains to be explained.
The approximate formula  \eqn{108}  gives no information on the
precise spectral nature of  the  ``small oscillations''
in $S_f(n)$.
Perhaps further information can be extracted from the
Fourier transform of the interpolating function defined in  \S7.

Earlier computations of 
J. B. Keiper \cite{Ke92} also apply to the Li
coefficients, as noted by K. Maslanka. In our notation
Figure 1 of his paper 
plots the function 
$$
\bar{\lambda}_n := \lambda_{-n} 
- \left( \frac{1}{2} n \log n + 
\frac{1}{2}(\gamma - 1 - \log (2\pi))n  \right).
$$
over the range $1 \le n \le 7000$.
In view of  Lemma~\ref{le42} and  
Theorem~\ref{th51} below (for  $\pi_{triv}$), the functions
$\bar{\lambda}_n$ 
and $S_f(n)$ differ by a bounded quantity; in fact
the actual difference appears to be very small. 
One may therefore interpret Keiper's
Figure 1 as essentially picturing  $S_f(n)$ over this range;
it indicates that  $|S_f(n)| < 20$ holds for $1 \le n \le 7000.$

Other related work on Li coefficients includes that of F. Brown \cite{Br04},
who determined zero-free regions for 
Dirichlet and Artin  $L$-functions
expressed in terms of the sizes of generalized Li coefficients.
P. Freitas \cite{Fr05} relates zero-free regions to the positivity of 
a different generalization of Li coefficients.
While this paper was being completed, an asymptotic
formula for the Li coefficients $\lambda_n$ 
was announced by A. Voros~\cite{Vo04}, under the
Riemann hypothesis. Comparison of his formula with 
that obtained here in \eqn{532}
led to a simplification of the expression
for $C_1(\pi)$ in Theorem~\ref{th51}.

\paragraph{Acknowlegments.} 
The author thanks the referee for many helpful
comments and simplifications of some details in \S5.
He thanks K. Maslanka for communicating 
various numerical constants
and graphs of Li coefficients, and for pointing out the
work of J. Keiper. He thanks  
E. Bombieri, Xian-Jin Li  and A. M. Odlyzko for helpful comments,
and  M. Coffey for 
communicating preliminary versions of his work.
Much of  this work was done while the author was
at AT\&T Labs-Research, whom he thanks for support.
%
%
%
%

\section{Li Coefficients for Automorphic $L$-functions}
\hsp
We recall basic facts about principal $L$-functions
$L(s, \pi)$ attached to irreducible cuspidal unitary
automorphic representations of $GL(N)$,
as in Jacquet \cite{Ja79} and 
Rudnick and Sarnak\cite[Sect. 2]{RS96}, see also 
Cogdell \cite{Cog03} and Gelbart and Miller \cite[Sect. 7.2]{GM03}.
  These $L$-functions 
are associated to $GL(n, \QQ)\backslash GL(n, \AAA_{\QQ})$,
and they are more precisely denoted $L(s, \pi, \rho)$ in
which the Langlands
$L$-group ${}^{L}G =GL(N, \CC)$ and 
$\rho : {}^{L}G \to GL(N, \CC)$ is the
standard  representation. 

For the trivial representation $\pi_{triv}$ of $GL(1)$
we have the completed automorphic $L$-function 
$\Lambda(s, \pi_{triv}) = \pi^{-\frac{s}{2}} \Gamma(\frac{s}{2}) \zeta(s)$.
This function has simple poles at $s=0$ and $s=1$.
Aside from this representation, all other $\Lambda(s, \pi)$ are
entire functions. 

Each completed automorphic $L$-function 
$\Lambda(s, \pi)$ has an Euler product factorization
\beql{N201}
\Lambda(s, \pi) :=  Q(\pi)^{\frac{s}{2}}L_{\infty}(s, \pi) L(s, \pi).
\eeq
Here $Q(\pi)$ is a positive integer called the {\em conductor}
of the representation $\pi$, and the {\em archimedean factor} is 
\beql{N202}
L_{\infty}(s, \pi) =  \prod_{j=1}^N \Gamma_{\RR}( s +  \kappa_j(\pi)),
\eeq
in which  $\kappa_j(\pi)$ are certain constants and 
\beql{N202a}
\Gamma_{\RR}(s) := \pi^{-\frac{s}{2}}\Gamma(\frac{s}{2}).
\eeq
The {\em (finite) $L$-function} $L(s, \pi)$ is given by an Euler product
over the finite places  
\begin{eqnarray*}
L(s, \pi) &= &\prod_p 
\prod_{j=1}^N \left(1 - \alpha_{p, j}(\pi)p^{-s}\right)^{-1}. \\
&= &\sum_{n=1}^{\infty} a_n(\pi) n^{-s}.
\end{eqnarray*}
This Euler product and its associated
Dirichlet series   
converges absolutely in a half-plane specified below. 
The functions $\Lambda(s, \pi)$ satisfy 
a functional equation
\beql{N203a}
\Lambda(s, \pi) = 
\epsilon(\pi) \Lambda( 1- s, \pi^{\vee}),
\eeq 
in which   $\epsilon(\pi)$ is a constant of
absolute value one, 
and $\pi^{\vee}$ denotes the contragredient representation.
The contragredient representation has $L$-function
\beql{N204}
L(s, \pi^{\vee}) = \sum_{n=1}^{\infty} \overline{a_n(\pi)} n^{-s}
\eeq
and archimedean factor 
\beql{N204a}
L_{\infty}(s, \pi^{\vee}) = L_{\infty}(s, \pi),
\eeq
and has  conductor $Q({\pi}^\vee) = Q(\pi)$.
The functional equation implies
 that $\epsilon(\pi^{\vee}) = \overline{\epsilon(\pi)}$.
The functions $\Lambda(s, \pi)$  are bounded in vertical strips, 
with exponential
decay as $|\Im(s)| \to \infty$, with $-B < \Re(s) < B$ for any
fixed $B$.

We define the $\xi$-function $\xi(s, \pi)$ associated to $\pi$ by 
\beql{N205}
\xi(s, \pi) := s^{-e(0, \pi)}(s-1)^{-e(1, \pi)}
\left(\frac{1}{\sqrt{ (-1)^{e(\frac{1}{2}, \pi)}\epsilon(\pi)}}
\Lambda(s, \pi) \right)~, 
\eeq
where $e(s_0, \pi)$ denotes the order of a zero or pole  of $\Lambda(s, \pi)$
at $s=s_0$, with poles being assigned negative orders.
 We have $e(0, \pi)= e(1, \pi)$ by the functional equation,
and this definition ensures that $\xi(s, \pi)$ is holomorphic and
nonzero at $s= 0$ and $1$. In this definition the square roots
must be chosen consistently so that 
$$
\sqrt{ (-1)^{e(\frac{1}{2}, \pi)}\epsilon(\pi) }~ \cdot
\sqrt{(-1)^{e(\frac{1}{2}, \pi)}\epsilon(\pi^{\vee})} = 1.
$$
There remains a choice of sign,  which can be  removed by 
the requirement that 
$\xi(\frac{1}{2}+ it, \pi)>0$ hold for small positive
$t$, as justified in the result below. 
  For the trivial representation
$\pi_{triv}$ on $GL(1)$ we have 
$e(0, \pi_{triv})=e(1, \pi_{triv}) = -1$,
and $\xi(s, \pi_{triv}) = 2 \xi(s)$.
This  convention is forced
if  we wish to have entire functions in all cases, for we must
remove the poles 
at $s=0$ and $s=1$ for the case $\pi_{triv}$.

The following theorem  collects together analytic facts about
automorphic $L$- functions.

%
%

\begin{theorem}~\label{th21}
Let $\pi$ be a   irreducible cuspidal
unitary automorphic representation of $GL(n)$ over $\QQ$.

(1) The ordinary Dirichlet series 
$L(s, \pi)= \sum_{n=1}^\infty a_n(\pi) n^{-s}$ converges 
absolutely in the half-plane $\Re(s) > 1$.
We have for all $n \ge 1$ that 
\beql{N206}
|a_n(\pi)| \le C(\pi) d(n) n^{\frac{N}{2}}, 
\eeq
for some $C(\pi) > 0$ depending on $\pi$, and $d(n)$
is the number of divisors of $n$. 

(2) The archimedean factors  $\Gamma_{\RR} (s + \kappa_j(\pi))$
in the Euler product $\Lambda(s, \pi)$ satisfy
\beql{N206a}
\Re( \kappa_{j}(\pi)) > -\frac{1}{2}.
\eeq
The quantities $\{\kappa_{j}(\pi): ~ 1\le j \le N\}$ are permuted under
complex conjugation, 
so that
\beql{N203}
L_{\infty}(s, \pi) = \overline{ L_{\infty}(\bar{s}, \pi)}.
\eeq

(3) The zeros of $\Lambda(s, \pi)$ all lie in the 
open critical strip $0 < \Re(s) < 1$. In particular
$\Lambda(s, \pi)$ is non-vanishing  on the lines
 $\Re(s)=0$ and $\Re(s)=1$. 

(4) The counting function  $N_{\pi}^{+}(T)$ (resp. $N_{\pi}^{-}(T)$) for
zeros of $\Lambda(s, \pi)$ with $0 \le \Im(\rho_{\pi}) < T$
(resp.
$-T \le \Im(\rho_{\pi}) \le 0$) each satisfy
\beql{N208}
N_{\pi}^{\pm}(T) = \frac{N}{2\pi} T \log T + \frac{1}{2}C_0(\pi) T + 
O ( \log T)
\eeq
as $T \to \infty$. 
Here 
\beql{N208a}
C_0(\pi) =\frac{1}{\pi}\log Q(\pi) - \frac{N}{\pi}(1+ \log (2\pi)),
\eeq  
in which  $Q(\pi)$ is the conductor of $\pi$, and the
$O$-constant depends on $\pi$.

(5) $\xi(s, \pi)$ satisfies the   functional equation 
\beql{N207}
\xi(s, \pi) = (-1)^k \xi(1-s, \pi^{\vee}),
\eeq
where $k= e(\frac{1}{2}, \pi) = e(\frac{1}{2}, \pi^{\vee})$ 
is the order of the zero of  $\xi(s, \pi)$ 
at $s=\frac{1}{2}$.
It is real-valued  on the critical line $\Re(s) = \frac{1}{2},$
so the multiset $Z(\pi)$ is invariant under the map
$\rho \mapsto 1 - \bar{\rho}.$

(6) The function $\xi(s, \pi)$ is an entire function of
order one and maximal type. It is bounded in vertical strips
$-B < \Re(s) < B$ for any finite $B$, and has rapid decrease
there as $|\Im(s)| \to \infty$. 
\end{theorem}

\paragraph{Proof.}
(1) Those irreducible cuspidal automorphic
representations $\pi$ that arise as  subrepresentations
of the right regular representation on 
$L^2(GL(n, \QQ) \backslash GL(n, \AAA_{\QQ}),
d\mu_{Haar})$ are necessarily unitary. 
We use the bound of 
Jacquet and Shalika \cite[Theorem 5.3]{JS81},
that Rankin-Selberg convolutions $L(s, \pi_1 \times \pi_2)$
of unitary irreducible cuspidal 
automorphic representations $\pi_1, \pi_2$ have Dirichlet
series that converge in $\Re(s) > 1$. We take $\pi_1=\pi$
and $\pi_2$ to be the trivial representation.
(See also \cite[Sect. 5.11, 5.12]{IK04}.)

Jacquet and Shalika \cite[Corollary 2.5]{JS81} also derive the bound  
that all unramified primes in the Euler
product have  $|\alpha_{j, p}| < \sqrt{p}$. 
All but finitely many primes are unramified, and therefore 
multiplying out the Euler product yields the
bound \eqn{N206}.
Better bounds are known for $|\alpha_{j, p}|$ which
give larger regions of absolute convergence.
Rudnick and Sarnak \cite[Prop. 5.1]{RS96} show
for irreducible cuspidal automorphic repesentations
$\pi$ of $GL(N)$ unramified at $p$ satisfy
$$
|\alpha_{j, p}(\pi)| \le p^{\frac{1}{2}- \frac{1}{N^2 + 1}}.
$$
The {\em generalized Ramanujan conjecture} 
formulated in  \cite{IS00} asserts that each
$|\alpha_{j, p}| = 1$ at unramified places $p$. It is
known to be true for principal 
$L$-functions over $GL(1)$ (Dirichlet $L$-functions).

(2) The inequality $\Re(s) > - \frac{1}{2}$ is
established in 
Rudnick and Sarnak ~\cite[eqn. (2.5) and Sect. 5.3]{RS96}.
There exist
ramified archimedean representations 
over $GL(2)$ for which some $\Re(\kappa_{j}(\pi)$ is
arbitrarily large. 
However
Luo, Rudnick and Sarnak \cite{LRS99} show that
unramified archimedean representations
\footnote{A representation is {\em unramified} if it is identically
one on the maximal compact subgroup of all the archimedean
components of $GL(N)$. Note that a  Dirichlet character $\chi$
for $GL(1)$ with $\chi(-1)=-1$ corresponds to an
automorphic  representation for $GL(1)$ that is 
ramified at the archimedean
place, and has  $\kappa_1(\chi)=1$.}
satisfy
the stronger bound
$$
|\Re(\kappa_j(\pi))| \le \frac{1}{2} - \frac{1}{N^2 + 1}.
$$
The Ramanujan conjecture at the archimedean places
asserts for an unramified representation that all 
$\Re( \kappa_j(\pi)) =0$. Its truth would imply that
the bound \eqn{N206a} could be improved to
$\Re(\kappa_j(\pi)) \ge 0$ in the general case.

The symmetry under complex conjugation of the
$\kappa_j(\pi)$ holds because the local factors
in the Euler product are equivalent to  
those of the contragredient
(see \cite[Theorem 2]{GK74} and \cite[Sect. 2.2]{RS96}),
which gives \eqn{N203}.

(3)  Iwaniec and Kowalski \cite[Theorem 5.42]{IK04}
show non-vanishing of 
cuspidal automorphic $L$-functions for $GL(N)$
on the line $\Re(s) =1$, and also obtain a
zero-free region inside the critical strip.
The nonvanishing on $\Re(s) =0$ comes from the functional equation
\eqn{N203a} for $\Lambda(s, \pi)$.

(4) The asymptotic formula for $N_{\pi}(T)$ is
essentially determined by the archimedean factors in the Euler
product for $\Lambda(s, \pi).$ The counting result
the zero density from $-T$ to $T$ appears
as \cite[Theorem 5.8]{IK04}, in which
$Q= Q(\pi)$ is the conductor of $\pi$ (see \cite[Sec. 5.1]{IK04}).
The error term given is  $O\left( \log \fq(\pi, iT)\right)$,
with an absolute constant, which involves 
the {\em analytic conductor}
\beql{214a}
\fq(\pi, s) := Q(\pi) \prod_{j=1}^N ( |s + \kappa_j(\pi)| + 3).
\eeq
Since we regard $\pi$ as fixed, this yields
$O(\log T)$, with the $O$-constant depending on $\pi$.
To get the bound (4) on the upper and lower critical
strip separately, the contour
integral proof in \cite{IK04} must be modified to split  into two
contours with a cut along the real axis, which
goes off it in small circles 
to avoid poles at the trivial zeros. 
(We omit the details.)

(5) The definition of $\xi(s)$ ensures that the
functional equation \eqn{N206} holds up to a sign,
which depends on the multiplicity $(\bmod ~2)$ of
a zero of $\Lambda(s, \pi)$ at $s= \frac{1}{2}$.
Now the  functions $\Lambda(s, \pi)$ have the symmetry
\beql{N209}
\Lambda (s, \pi) = \overline{ \Lambda(\bar{s}, \pi^{\vee} ) }.
\eeq
which follows from \eqn{N203} and \eqn{N204}. 
In consequence \eqn{N205} gives
$$
\xi(s, \pi) = (-1)^{e(\frac{1}{2}, \pi)} 
\overline{ \xi(\bar{s}, \pi^{\vee}) }.
$$
We deduce that
$$
 \xi( \frac{1}{2} + it, \pi) = 
(-1)^{e(\frac{1}{2}, \pi)}  \xi ( \frac{1}{2} - it, \pi^{\vee}) =
\overline{ \xi ( \frac{1}{2} + it, \pi) }.
$$
Thus $\xi(s, \pi)$ is real on the critical line.
The invariance of the multiset $Z(\pi)$ under $\rho \mapsto 1 - \bar{\rho}$
now follows from the reflection principle.

(6) Godement  and Jacquet \cite[Theorem 13.8]{GJ72} showed for
cuspidal automorphic representations $\pi$
that $\Lambda(s, \pi)$
is meromorphic and  bounded in vertical strips, with
finitely many poles, and that 
the ordinary Dirichlet series 
$L(s, \pi)$ has a nonempty half-plane of absolute convergence,
(See also  Jacquet \cite[Theorem 6.2]{Ja79}). 
The assumption that $\pi$ is irreducible and cuspical 
gives that they are entire functions, aside from the 
trivial representation $\pi_{triv}||\cdot||^s$ of $GL(1)$,
whose singularities are simple poles at $s=0, 1$.
It follows that $\xi(s, \pi)$ is an entire function in all cases.

We assert that in the vertical strip $-B \le \Im(s) \le B$ 
the function $\xi(s, \pi)$ satisfies 
$$
|\xi(s, \pi)| \le C(\pi, B) e^{- N \frac{\pi}{2}|\Im(s)|} 
$$
as $|\Im(s)| \to \infty$. To see this,
observe that the exponential decay on vertical lines
holds for $\Re(s) \ge 1 + \epsilon$, coming from exponential
decay of the archimedean factor and the bound $|L(s, \pi)| = O(1)$.
It then holds for $\Re(s) \le  - \epsilon$
by the functional equation, and then in between by the
Phragm\'{e}n-Lindel\"{o}f principle.

In the half-plane $\Re(s) \ge 1 + \epsilon$ we have
$|\Lambda(s, \pi)|= O(1)$ because its Dirichlet series
representation converges absolutely. It follows that
on this half-plane 
$|\Lambda(s, \pi)| = O (e^{ N |s| (\log |s|+1)})$,
with the growth rate coming from the archimedean  factors,
which can be bounded by Stirling's formula.
The functional equation now shows that the same bound
holds on the half plane $\Re(s) \le -\epsilon$. Since 
 $|\Lambda(s, \pi)|$ is known to be bounded on the vertical strip
$ -\epsilon \le \Re(s) \le 1+\epsilon$, 
it follows that $|\Lambda(s)| = O( e^{N |s| (\log (|s| +2)})$
for all $s$. Thus $\Lambda(s)$ is
an entire function of order one.
It is of maximal type because its growth-rate is 
faster than $e^{R|s|}$ for any finite $R$ along
the positive real axis, coming from the archimedean factors.
The function  $|\xi(s, \pi)|$ inherits these properties,
since it differs by at most
a polynomial factor from $\Lambda(s, \pi)$.
$~~~\bsq$

We let $Z(\pi)$ denote the multi-set of zeros of $\xi(s, \pi)$
(counted with multiplicity)
which is the same as that of $\Lambda(s, \pi)$ except
possibly at $s=0$ and $s=1$. 

\begin{lemma}~\label{Nle21}
For any principal  $L$-function $L(s, \pi)$ for $GL(N)$
the power sums 
\beql{N210}
\sigma_n(\pi) := \sum_{\rho \in Z(\pi)}{}^{'} \frac{1}{\rho^n},~~~ n \ge 1, 
\eeq
are  absolutely convergent
for $n \ge 2$,  and are 
$*$-convergent for $n=1$. The real parts of these sums
are absolutely convergent for all $n \ge 1$.
\end{lemma}

\paragraph{Proof.}
The absolute convergence for $n \ge 2$ follows easily
from the zero-counting bound \eqn{N208}.

The $*$-convergence for $\sigma_1(\pi)$ follows from the asymptotics
\eqn{N208} for  the zeros. Here the  zeros below and above the
real axis are paired in increasing order of their imaginary
parts (in absolute values) and using the asymptotic
formulas in Theorem~\ref{th21}(4) we get $*$-convergence
by partial summation. It is important that the
remainder term in \eqn{N208} be $O(T^{1 -\delta})$ for
some $\delta > 0.$
(In general the  zeros are
not symmetric about the real axis. The symmetry
$\rho \mapsto 1 - \bar{\rho}$ is of no help in proving
$*$-convergence.)
Finally, for any zero $\rho= \beta + i \gamma$ we have 
$$
\Re \left(\frac{1}{\rho}\right)= \frac{\beta}{\beta^2 + 
\gamma^2} = O \left(\frac{1}{\gamma^2}\right),
$$
which  gives absolute convergence for $n=1$ of the 
real parts of the terms in the sum \eqn{N210}.
$~~~\bsq$

\begin{lemma}~\label{Nle22}
For all irreducible cuspidal unitary automorphic representations
on $GL(N)$ the sums 
\beql{N211}
\lambda_n(\pi):= \sum_{\rho \in Z(\pi)} 1 - 
\left( 1 - \frac{1}{\rho}\right)^n
\eeq
are $*$-convergent for all $n \in \ZZ$. They are given by 
\beql{N212}
\lambda_n(\pi) = \sum_{j=1}^n (-1)^{j-1} {\binom{n}{j}} \sigma_{j}(\pi),
\eeq
with  $\lambda_0(\pi)  =0$. They satisfy
\beql{N213}
\lambda_{-n}(\pi) = \overline{\lambda_n(\pi)} = \lambda_n(\pi^{\vee}).
\eeq
\end{lemma}

\paragraph{Proof.}
We have
$$
\sum_{|\rho| < T} 1 - \left( 1 - \frac{1}{\rho}\right)^n =
\sum_{|\rho| < T} \frac{n}{\rho} + 
\sum_{j=2}^n (-1)^{j-1}{\binom{n}{j}} 
\left(  \sum_{|\rho| < T}\frac{1}{\rho^j} \right).
$$
On letting  $T \to \infty$, Lemma~\ref{Nle21} shows
that the first sum on the right $*$-converges and the 
second sum on the right converges absolutely.
This gives \eqn{N212}.

The symmetry $\lambda_{-n}(\pi) = \overline{\lambda_n (\pi)}$
is inherited from the symmetry that if $\rho= \beta + i\gamma$
is a zero of $\xi(s)$, then so is $1- \bar{\rho} = 1-\beta +i\gamma$,
with the same multiplicity. 
This holds by the reflection principle since $\xi(s)$ is
real on the critical line.
That is, the multiset of
zeros is invariant under the map $\rho \mapsto 1 - \bar{\rho}.$
Now we have 
$$
1- (1 - \frac{1}{\rho} )^{-n} = 1- (\frac{\rho -1}{\rho})^{-n}
= 1- (\frac{-\rho}{1-\rho})^n = 
1- \overline{ \left(1-  \frac{1}{1-\bar{\rho}}\right)^n }.
$$
and this gives $\lambda_{-n}(\pi) = \overline{\lambda_n(\pi)}$. 

For the contragredient representation, we have the symmetry
$\Lambda (s, \pi) =\overline{\Lambda(\overline{s}, \pi^\vee)}.$
in \eqn{N209}. This shows that the zero sets of the two
completed $L$-functions are complex conjugate, i.e.
$Z(\pi^{\vee}) = \overline{Z(\pi)}.$ This then yields the 
other relation
$\lambda_n(\pi^{\vee}) = \overline{\lambda_n(\pi)}$.
$~~~\bsq$

We can now state a general version of Li's criterion.

\begin{theorem}~\label{Nth22}
Let $\pi$ be an irreducible cuspidal
unitary automorphic representation of $GL(n)$. The following
conditions are each equivalent to the 
Riemann hypothesis for $\xi(s, \pi)$.


(1) For all $n \ge 1$, 
\beql{N217a}
\Re \left( \lambda_n(\pi) \right) \ge 0.
\eeq

(2) For each $\epsilon > 0$, there is a positive
constant $C(\epsilon)$ such that 
\beql{217b}
\Re \left( \lambda_n(\pi) \right) \ge - C(\epsilon) e^{ \epsilon n}
~~~\mbox{for all}~~~ n \ge 1.
\eeq

(3) The generalized Li coefficients $\lambda_n(\pi)$ satisfy
\beql{N214}
\lim_{n \to \infty} |\lambda_n(\pi)|^{\frac{1}{n}} \le 1.
\eeq
\end{theorem}

\paragraph{Proof.}
Theorem~\ref{th21} gives that 
the multiset $Z(\pi)$ omits the values $0$
and $1$ and is invariant under the symmetry $\rho \mapsto 1 - \bar{\rho}.$
The equivalence of conditions (1) and  (2) to the
Riemann hypothesis for $Z(\pi)$  follows from the Corollary
in Theorem 1 of \cite{BL99}. 

It remains to show  the equivalence of (3)
to the Riemann hypothesis.  The $\lambda_n(\pi)$ are
identifiable as the power series coefficients around 
$z=0$ of
\beql{N215}
\frac{d}{dz}  \log \xi(\frac{-z}{1-z}, \pi) =
\frac{-1}{(1-z)^2} \frac{\xi'}{\xi} (\frac{-z}{1-z}, \pi) 
= \sum_{n=0}^{\infty}
\lambda_{n+1}(\pi) z^n.
\eeq
The map $z \mapsto s = -\frac{z}{1-z}$ conformally maps the 
unit disk to the half-plane
$\Re(s) \ge \frac{1}{2}$, and the point
$z=0$ corresponds to $s=1$, a point
 where $\frac{\xi'}{\xi}(s)$
is holomorphic.  Assuming the Riemann hypothesis holds
for $\xi(s, \pi)$, the function in \eqn{N215} is holomorphic
in the unit disk, hence its power series coefficients
satisfy \eqn{N214}. Conversely, if \eqn{N214} holds, then
the function \eqn{N215} is holomorphic in the unit disk,
hence $\frac{\xi'}{\xi}(s, \pi)$ has no singularity
in $\Re(s) > \frac{1}{2}$. The functional equation then shows
it has no singularity in $\Re(s) < \frac{1}{2}$ so the
Riemann hypothesis holds for $\xi(s, \pi).$  $~~~\bsq$

%
%
%
\section{Li Coefficients and Weil's Quadratic Functional}
\hsp
A. Weil \cite{We52}  
formulated the ``explicit formula'' of prime number
theory in terms of  distributions, and
using the Fourier transform on the real line.  He gave a necessary and
sufficient condition for  the Riemann
hypothesis in terms of the positivity of a quadratic
functional on a suitable space of test functions on
the real line, contained in the Hilbert space $L^{2}(\RR, du)$. 
There is a natural
generalization of Weil's quadratic functional associated
to the completed $L$-function $\Lambda(s, \pi)$  
of any irreducible cuspidal  automorphic representation $\pi$
of $GL(N)$, which we shall conisder here.

We first note that, using an
exponential change of variable ($x=e^{u}$), 
the  ``explicit formula'' 
is expressible in 
terms  of   test functions 
contained in the Hilbert space $L^2( \RR_{>0}, \frac{dx}{x})$
on the positive real line $\RR_{>0}$ in the $x$-variable,
with the Fourier transform replaced by the Mellin transform.
This is the framework for the ``explicit
formula''  taken in \cite{Bo99}, \cite{BL99}.
Now we make  a second change of variable, using the
Mellin transform 
$$\hat{f}(s)= \int_{0}^{\infty} f(x) x^{s}\frac{dx}{x},$$
to transfer Weil's functional to a quadratic functional on
a space of test functions contained inside the Hilbert space
$L^2(\frac{1}{2} + i\RR, \frac{dt}{2\pi})$; these test
functions are the Mellin transforms  on the critical line 
of the test functions
on  $\RR_{>0}$ above. (The Mellin transform extends to an
isometry between $L^{2}(\RR_{>0}, \frac{dx}{x})$
and $L^{2}(\frac{1}{2}+i\RR, \frac{dt}{2\pi}).$)
We treat the ``explicit formula'' in these coordinates
(a viewpoint taken in  Burnol~\cite{Bu98}).
The resulting  test functions in the $s$-variable
 have the property of being
analytic in some open domain that includes the 
critical line $\Re(s) = \frac{1}{2}$ in its interior.
If the Riemann hypothesis is not assumed, the set
of  test functions must be further restricted to 
functions  analytic in a region containing the closed  critical strip. 

We consider as test functions the vector space 
$\sA$ of all 
functions $F(s)$ holomorphic
in the strip $0 < \Re(s) < 1$ which satisfy a uniform  growth bound
$F(s) = O( \frac{1}{|s|})$ in the strip outside
\footnote{The region $|\Im(s)| \le 1$ is omitted
here  to  avoid the points $s=0$ and $s=1$, because every 
nonzero element in the Li class $\sL$ of test
functions defined later necessarily has a pole at
one or both of these points.}
the region 
$|\Im(s)| \le 1$, with $O$-constant depending on
the function.  The class $\sA$ is closed under 
the action of the involution $\tilde{G}(s) := G(1-s)$.
The functions in $\sA$ are completely determined by their
values on the critical line $\Re(s) = \frac{1}{2}$ by analytic
continuation, and the growth bound ensures that 
they belong to $L^2(\frac{1}{2} + i\RR, \frac{dt}{2\pi})$.
Under the inverse Mellin transform they convert to 
a class of test functions contained in 
$L^2(\RR_{>0} , \frac{dx}{x})$ of the type considered in \cite{BL99};
these functions are smooth away from the point $x=1$.

Given $F(s), G(s) \in \sA$ we define the {\em
Weil scalar product}  associated to
the automorphic representation $\pi$ by 
\beql{301}
\langle F, G \rangle_{\sW(\pi)} := \sum_{ \rho \in Z(\pi) }
F(\rho) \overline{G( 1 - \bar{\rho})}.
\eeq
The sum on the right counts zeros with
multiplicity, and it converges absolutely due to the 
growth bound on $F$ and $G$ for large $|s|$.
This scalar product is linear in the first factor and
conjugate-linear in the second factor.  The
multiset of automorphic $L$-function zeros $Z(\pi)$
(counting multiplicities)  is invariant under
the involution $\rho \to 1 - \bar{\rho}$.
This yields 
the Hermitian symmetry  
$$
\langle F, G \rangle_{\sW(\pi)} = 
\overline{ \langle G, F \rangle}_{\sW(\pi)} .
$$
In the appendix we clarify the relation between this
definition and Weil's definition of the 
quadratic functional. 

The Riemann hypothesis for $\Lambda(s, \pi)$ implies
that the Weil scalar product is positive semidefinite on
the test function vector space $\sA$. To see
this we note that  
$\rho= 1 - \bar{\rho}$ holds if and only if 
$\rho$ lies on the critical line $\Re(s) = \frac{1}{2}.$
The Riemann hypothesis implies that for all $F(s) \in \sA$, 
$$
\langle F, F \rangle_{\sW(\pi)} = \sum_{ \rho \in Z(\pi) } 
F(\rho) \overline{F(1 - \bar{\rho})}
= \sum_{ \rho \in Z(\pi) }
F(\rho) \overline{F(\rho)} = \sum_{ \rho \in Z(\pi) }
|F(\rho)|^2 \ge 0.
$$
The  Weil scalar product is not positive definite on
the full class $\sA$ since $F(s) := \xi(s, \pi) \in  \sA$ with 
$\langle F, F \rangle_{\sW(\pi)}=0$. In the converse
direction,  Weil showed for various $L$-functions
that the positive semidefiniteness of the Weil scalar product
on suitable subsets of the test functions in $\sA$.
implies the Riemann Hypothesis for $\Lambda(s, \pi)$.
A number of criteria of this sort are known,
using different collections of test functions; essentially
one needs the test function set to be sufficiently large to 
separate all the zeros $\rho = \rho_{\pi}$
(counting  multiplicities).

 We define the {\em Li class} $\sL$ of test functions to be 
the set of rational functions in  the function field $\CC(s)$
that vanish 
at infinity (on the Riemann sphere) and whose polar
divisor is contained in the set $\{ 0, ~ 1\}$. The class
$\sL$ is closed under addition, multiplication and
scalar multiplication, but does not contain the constant
functions.
A vector space basis for $\sL$ consists of $\{\frac{1}{s^n}:~
n \ge 1\}$ and $\{\frac{1}{(1-s)^n}:~ n \ge 1\}$. 

The special test functions $G_n(s) \in \sL$
corresponding to the Li coefficients  are
\beql{303}
G_n(s) := 1 - (1 - \frac{1}{s})^n ~~~\mbox{for}~~~ n \in \ZZ.
\eeq
The set of all $G_n(s)$, excluding $G_0(s) \equiv 0$,
forms a vector space basis of $\sL$. Indeed  the 
change-of-basis matrix $G_n(s)$ $(n \ge 1)$
relating these functions to $\{\frac{1}{s^n}:~ n \ge 1\}$
is an upper triangular unipotent matrix, and similarly
for $G_{-n}(s)$ $(n \ge 1)$ and $\{\frac{1}{(1-s)^n}:~ n \ge 1\}$.
The class $\sL$ is contained in $\sA$, because the vanishing
condition at $\infty$ implies a bound $F(s) = O(\frac{1}{|s|})$
uniformly in the region $|\Im(s)| \ge 1$.   
The class $\sL$  has 
the property that every nonzero member of it
has a pole at either $s=0$ or $s=1$, or both.

The following result computes the Weil scalar product with
respect to this basis.

%
%

\begin{theorem}~\label{th31}
Let $\pi$ be an irreducible cuspidal unitary automorphic representation
of $GL(N)$, with  associated $\xi$-function
$\xi(s, \pi)$
and with Weil scalar product $\langle \cdot, \cdot \rangle_{\sW(\pi)}$.
For the Li test functions  $G_n(s) = 1 - (1 - \frac{1}{s})^n$
there holds 
\beql{304}
\langle G_n, G_m \rangle_{\sW(\pi)} = 
\lambda_n(\pi) + \lambda_{-m}(\pi) - \lambda_{n-m}(\pi).
\eeq
In particular 
\beql{304a}
||G_n||_{\sW(\pi)}^2 = \lambda_n(\pi) + \lambda_{-n}(\pi) 
= 2 \Re(\lambda_n(\pi)).
\eeq
\end{theorem}

\paragraph{Proof.}
By definition
$$
\lambda_n(\pi) = \sum_{\rho \in Z(\pi)}{}^{'} G_n(\rho).
$$ 
Since each function
$G_n(s)$ is real on the real axis, the reflection principle
gives
$$
\overline{ G_n(1 - \bar{\rho}) } = G_n(1 - \rho).
$$
We also have the identity 
\beql{205}
G_{-m}(s) = 1 - (1 - \frac{1}{s})^{-m} 
= 1 - ( 1- \frac{1}{s-1})^m = G_m(1-s)
\eeq 
and  
\beql{206}
G_n(s) G_{-m}(s) = G_n(s) + G_{-m}(s) - G_{n-m}(s).
\eeq
Combining all these gives
\begin{eqnarray*}
\langle G_n, G_m \rangle_{\sW(\pi)} &= & 
\sum_{\rho \in Z{\pi}} G_n(\rho) G_m(1 - \rho) 
= \sum_{\rho} G_n(\rho)G_m(\rho) \\
&= &\sum_{\rho \in Z(\pi)}{}^{'} ( G_n(\rho) + 
G_{-m}(\rho) - G_{n-m}(\rho)) \\
&= &\sum_{\rho \in Z(\pi)}{}^{'} G_n(\rho) + 
\sum_{\rho \in Z(\pi)}{}^{'} G_m(\rho) + 
\sum_{\rho \in Z(\pi)}{}^{'} G_m(\rho) \\
&= & \lambda_n(\pi) + \lambda_{-m}(\pi) - \lambda_{n-m}(\pi).
\end{eqnarray*}
the required result. $~~~\bsq$.

The  Li class $\sL$ is  large enough 
to characterize the Riemann hypothesis in terms of the
semidefiniteness of the Weil scalar product.
Theorem~\ref{th31} shows that
semidefiniteness implies that all $\Re(\lambda_n) \ge 0$.
The Riemann hypothesis then follows from  
 Theorem~\ref{Nth22}(2).

If the Riemann hypothesis holds,
then the  Weil scalar product is positive definite on $\sL$.
For it  $F(x) \in \sL$ is nonzero, then
since it is a rational function it has finitely many
zeros, 
so cannot vanish at all points of $Z(\pi).$
Thus 
$$
||F||_{\sW(\pi)}^2 = \sum_{\rho \in Z(\pi)} |F(\rho)|^2 > 0.
$$
The vector space $\sL$ is then a pre-Hilbert space with
this Hermitian scalar product, so it  can be completed
to a Hilbert space $\sH_{\sL}(\pi).$


%
%
%
\section{Arithmetic Formula for Li  Coefficients}
\hsp
We relate the Li coefficients to
sums of the form 
$\sum_{\rho \in Z(\pi)}^{'} f(\rho)$ for various
functions $f(s)$. In particular we split $\lambda_n(\pi)$
into two parts $S_{\infty}(n, \pi)$
and $S_{f}(n, \pi)$ given in Lemma~\ref{le42} below.

The Li coefficients are
expressible in terms of Laurent series
coefficients around the point  $s=1$ of
of the logarithmic derivative of the  $\xi$-function
$\xi(s, \pi)$.  These coefficients are given
in terms of power sums of zeros. 

\begin{lemma}~\label{le41}
The power series expansion around
$s=1$ of $\frac{\xi'}{\xi}(s, \pi)$ is 
\beql{N405} 
\frac{\xi'}{\xi}(s+1, \pi) = \sum_{n=0}^{\infty}
(-1)^{n} \sigma_{n+1}(\pi^{\vee}) s^n.
\eeq
in which the power sum $\sigma_{n}(\pi^{\vee}) = 
\sum_{\rho \in Z(\pi^{\vee})}^{'} \frac{1}{\rho^n}.$
\end{lemma}

\paragraph{Proof.}
Since $\xi(s, \pi)$ is
entire of order one we  have the Hadamard product expansion
$$
\xi(s, \pi) = e^{A(\pi)s + B(\pi)}
\prod_{\rho \in Z(\pi)}( 1 - \frac{s}{\rho})e^{\frac{s}{\rho}}.
$$
Taking logarithms yields,
\begin{eqnarray*}
\log \xi(s, \pi) &= & A(\pi)s + B(\pi) +
\sum_{\rho \in Z(\pi)} \left(\log (1 - \frac{s}{\rho})+\frac{s}{\rho}\right) \\
&= & A(\pi)s + B(\pi) -
\sum_{\rho \in Z(\pi)} 
\left(\sum_{n=2}^\infty \frac{1}{n}\frac{s^n}{\rho^n} \right) \\
&= &  A(\pi)s + B(\pi) - \sum_{n=2}^\infty  \sigma_n(\pi)\frac{s^n}{n}.
\end{eqnarray*}
One can also deduce by a $*$-convergent rearrangement that  
$$
A(\pi) = - \sum_{\rho \in Z(\pi)}{}^{'} \frac{1}{\rho} = -\sigma_1(\pi). 
$$
Differentiating yields
$$
-\frac{\xi'}{\xi}(s, \pi) = \sum_{n=0}^\infty \sigma_{n+1}(\pi) s^n.
$$
The functional equation \eqn{N207}gives 
$\xi'(s, \pi) = (-1)^{k+1} \xi'(1-s, \pi^{\vee})$, 
from which we obtain 
$$
\frac{\xi'}{\xi}(s+1, \pi)= - \frac{\xi'}{\xi}(-s, \pi^{\vee}) 
= \sum_{n=0}^\infty (-1)^{n}
\sigma_{n+1}(\pi^{\vee}) s^n,
$$
as required. $~~~\bsq$.

We now write
\beql{N405b}
\frac{\xi'}{\xi}(s+1, \pi) = \left( \frac{1}{2} \log Q(\pi) 
+ \sum_{j=1}^N \frac{\Gamma_{\RR}^{'}}{\Gamma_{\RR}}( s +1+ \kappa_j(\pi))
\right)+ \frac{L'}{L}(s+1,\pi) - \frac{e(0, \pi)}{s+1} -\frac{e(1, \pi)}{s}
\eeq
We define Laurent series coefficients around $s=1$ of the first two
terms on the right, as follows.  

\begin{defi}~\label{de41}
{\em 
 The  coefficents
$\{\tau_n(\pi):~ n \ge 1\}$  are defined by
\beql{N405c}
 \sum_{n=0}^\infty \tau_{n}(\pi) s^j:=     \frac{1}{2} \log Q(\pi)   + 
\sum_{j=1}^N \frac{\Gamma_{\RR}^{'}}{\Gamma_{\RR}}( s +1+ \kappa_j(\pi))
\eeq    
in which we have using \eqn{N202a} that }
\beql{N405e}
\frac{\Gamma_{\RR}^{'}}{\Gamma_{\RR}}( s) = -\frac{1}{2}\log \pi
 + \frac{1}{2}\frac{\Gamma'}{\Gamma}(\frac{s}{2}). 
\eeq
\end{defi}

\begin{defi}~\label{de42}
{\em The  coefficients  $\{\eta_n(\pi):~n \ge 1\}$ are defined by 
\beql{N405d}
 \sum_{n=0}^\infty \eta_{n}(\pi) s^n :=
-\frac{L'}{L}(s+1,\pi) -\frac{\delta(\pi)}{s}
\eeq
where  $\delta(\pi) = 1$ if $\pi= \pi_{triv}$ and $\delta(\pi) = 0$
otherwise. (The minus  sign is present to agree with the convention
in \cite[eqn. (4.3)]{BL99}.)}
\end{defi}

We note that the coefficients
$\tau_n(\pi)$ are real-valued, as a consequence of Theorem ~\ref{th21}(2).
Also, the  coefficients $\eta_n(\pi)$ for 
 $\pi= \pi_{triv}$ are exactly the coefficients
$\eta_n$ given in \cite{BL99}.

%
%

\begin{lemma}~\label{le42}
Let $\pi$ be an irreducible cuspidal automorphic representation
on $GL(N)$ over $\QQ$. Then for all $n \ge 1$, 
\beql{N406a}
\lambda_n(\pi) = S_{\infty}(n,  \pi) - S_f(n, \pi^{\vee}) + 
\delta(\pi),
\eeq
in which
\beql{N406b}
S_{\infty}(n, \pi^{\vee}) = \sum_{j=1}^n {\binom{n}{j}} \tau_{j-1}(\pi^{\vee})
\eeq
and 
\beql{N406c}
 S_f(n, \pi^{\vee}) = \sum_{j=1}^n  {\binom{n}{j}} \eta_{j-1}(\pi^{\vee}),
\eeq
and $\delta(\pi)=1$  if $\pi= \pi_{triv}$ and $\delta(\pi)=0$    otherwise.
\end{lemma} 

\paragraph{Remark.} Note that $S_{\infty}(n, \pi)= S_{\infty}(n, \pi^{\vee})$
and $\delta(\pi)=\delta(\pi^{\vee})$.
Also, since $\lambda_{-n}(\pi) = \overline{\lambda_{n}(\pi)} = 
\lambda_n(\pi^{\vee}),$ we obtain from \eqn{N406a} applied to $\pi^{\vee}$
 that
\beql{406d}
\lambda_{-n}( \pi) = S_{\infty}(n, \pi) - S_f(n, \pi) + \delta(\pi).
\eeq

\paragraph{Proof.}
Comparison of Lemma~\ref{le41} for $\pi^\vee$
and the formula \eqn{N405b} yields
\beql{N407}
(-1)^{j} \sigma_{j+1}(\pi) = \tau_{j}(\pi^{\vee}) - \eta_{j}(\pi^{\vee})
+ (-1)^{j} \delta(\pi^\vee).
\eeq
Substituting this 
formula into \eqn{N212}, and using 
$-e(0, \pi) = \delta(\pi)=\delta(\pi^{\vee})$, we obtain
the desired result. 
$~~~\bsq$ \\

In the expression \eqn{N406a} the term $S_{\infty}(n,  \pi^{\vee})$
corresponds to the contribution of the archimedean primes
and $S_f(n, \pi^{\vee})$ corresponds to the finite primes. There are also
extra contributions from singularities at $s=0$ and $s=1$.
The $s=1$ contribution cancels against the singularity 
at $s=1$ of
$\frac{L'}{L}(s, \pi)$, but the $s=0$ singularity remains
and contributes the constant $-e(0, \pi)\delta(\pi^{\vee})$. 

For $\pi_{triv} = \pi_{triv}^{\vee}$ on $GL(1)$ the 
terms $S_{\infty}(n,  \pi)$ and $S_f(n, \pi)$
form parts of the arithmetic expression for the Li coefficients
$\lambda_n$ given in Theorem 2 of \cite{BL99}.
We have $\mu_1(\pi_{triv})=0$ and 
\beql{N407b}
S_{\infty}(n,  \pi_{triv}) = -\sum_{j=1}^n (-1)^{j+1}{\binom{n}{j}}
(1- \frac{1}{2^j}) \zeta^{\ast}(j),
\eeq
in which $\zeta^{*}(j)= \zeta(j)$ for $j \ge 2$,
and $\zeta^{*}(1)= \log(4\pi) + \gamma$, where $\gamma$
is Euler's constant. For the other term, writing
\beql{N402}
- \frac{\zeta'}{\zeta}(s+1) = 
\frac{1}{s} + \sum_{j=0}^{\infty} \eta_j s^j, 
\eeq
we have 
\beql{N407c}
S_f(n, \pi_{triv}) = \sum_{j=1}^n (-1)^{j-1}{\binom{n}{j}}\eta_{j-1}.
\eeq
In \cite{BL99} it is shown that
the  coefficients $\eta_j$ are given by
$$
\eta_j := \frac{(-1)^{j}}{j!}\lim_{N \to \infty} \left( \sum_{m=1}^N 
\frac{\Lambda(m) (\log ~m)^{j}}{m} 
- \frac{1}{j+1}(\log N)^{j+1} \right),
$$
in which  $\Lambda(m)$ is the von Mangoldt function. 

There exist formal arithmetic expressions
for  $S_{f}(n,  \pi)$ in the general case
that are similar in spririt to \eqn{N407c}.
Logarithmically differentiating
the Euler product for $L(s, \pi)$ in
the region  $\Re(s) > 1$ yields 
$$
-\frac{L'}{L}(s,\pi) = \sum_{n=1}^\infty \Lambda_{\pi}(n) n^{-s}
$$
in which for $n= p^m$ a prime power, 
\beql{N408}
\Lambda_{\pi}(p^m)= \frac{1}{m} \sum_{k=1}^N \left(\alpha_{k,p}(\pi)\right)^m,
\eeq
and $\Lambda_{\pi}(n)=0$ otherwise.
Now for all $\pi \ne \pi_{triv}$
using \eqn{N405d} (which is now analytic at $s=1$) we have  
$$
\eta_j(\pi) =- \frac{1}{j!}\frac{d^j}{ds^j}  
\left[\frac{L'}{L}(s,\pi) \right]|_{s=1}
$$
Using \eqn{N408} we obtain, formally,
\beql{N409}
\eta_j(\pi) \simeq 
\frac{(-1)^{j}}{j!} \sum_{m=1}^{\infty} 
\frac{\Lambda_{\pi}(m)(\log m)^{j}}{m}
\eeq 
The sum on the right 
at best converges conditionally ,
when viewed as $\lim_{T \to \infty} \sum_{n \le T}$.
The conditional convergence  is known to hold for 
$L$-functions on $GL(1)$. One can show 
$$
\eta_j(\pi) = \lim_{s \searrow 1}\frac{(-1)^{j}}{j!} \sum_{m=1}^{\infty} 
\frac{\Lambda_{\pi}(m)(\log m)^{j}}{m^s}
$$ 
with $s  \searrow 1$ along the real axis.
Explicit expressions for $\eta_j(\pi)$ for
various $L$-functions and modular forms are derived in
X. Li \cite{Li03}, \cite{Li04a}, \cite{Li04b}.

There exist simplified formulas for the coefficients 
$\tau_n(\pi)$ analogous to \eqn{N407b},
which involve values of the Hurwitz zeta function
$\zeta(s, z) = \sum_{j=0}^{\infty} \frac{1}{(n+z)^s}$
with $\Re(z) > 0$.

%
%

\begin{lemma}~\label{le43}
For an irreducible cuspidal automorphic representation $\pi$ on
$GL(N)$  and $n \ge 1$, there holds
\beql{N417}
\tau_n(\pi) = (-\frac{1}{2})^{n+1} \sum_{j=1}^n 
\zeta(n+1, \frac{\kappa_j(\pi)+1}{2}).
\eeq
In addition for $n=0$ there holds 
\beql{N418}
\tau_0(\pi) = \frac{1}{2} \log Q(\pi) - \frac{N}{2} \log \pi
+\frac{1}{2}\sum_{j=1}^N \frac{\Gamma'}{\Gamma}
\left(\frac{\kappa_j(\pi) - 1}{2}\right).
\eeq

where $\gamma$ is Euler's constant. 
\end{lemma}

\paragraph{Proof.}
We recall that $\Gamma_{\RR}(s) = \pi^{-\frac{s}{2}} \Gamma(\frac{s}{2})$,
so that 
\beql{N419}
\frac{\Gamma_{\RR}^{'}}{\Gamma_{\RR}}( s) = -\frac{1}{2}\log \pi
 + \frac{1}{2}\frac{\Gamma'}{\Gamma}\left(\frac{s}{2}\right).
\eeq
The digamma function $\psi(s) = \frac{\Gamma'}{\Gamma}(s)$
has the partial fraction expansion
\beql{N420}
 \psi(s+1)= \frac{\Gamma'}{\Gamma}(s+1)= - \gamma + \sum_{m=1}^\infty
\left( \frac{1}{m} - \frac{1}{s+m}\right),
\eeq
in which $\gamma \approx  0.5771$ is Euler's constant.
For $\Re( z) > -1$ 
we  define the power series coefficients $\{ \tilde{\tau}_n(z): n \ge 0\}$
by
\beql{N421}
  \sum_{n=0}^{\infty} \tilde{\tau}_n(z) s^n :=   
\frac{\Gamma'}{\Gamma}(s+1+z). 
\eeq
Here 
$$ 
\tilde{\tau}_n(z) = \frac{1}{n!} 
\frac{d^n}{ds^n}\left( \frac{\Gamma'}{\Gamma}(s+1+z)\right)|_{s=0}.
$$
We have 
\beql{N422}
\tilde{\tau}_0(z)= - \gamma + 
\sum_{m=1}^{\infty} \frac{z}{m(m+z)}.
\eeq
and, for $n \ge 1$, 
\beql{N423}
\tilde{\tau}_n(z)= (-1)^{n+1}
\sum_{m=1}^{\infty} \left(\frac{1}{m+z}\right)^{n+1},
\eeq
which follows 
by differentiating the partial fraction
expansion \eqn{N420} at $s= z+1$. This is essentially a
special value of the Hurwitz zeta function, for $n \ge 1,$ 
\beql{N424}
\tilde{\tau}_n(z) = (-1)^{n+1} \zeta(n+1, z+1).
\eeq
Comparing power series coefficients of the definition \eqn{N405c}
with those of \eqn{N419} yields, for $n \ge 1$,
\beql{N425}
\tau_n(\pi) = \frac{1}{2} \sum_{j=1}^N 
\tilde{\tau}_n( \frac{\kappa_j(\pi)-1}{2}) \left( \frac{1}{2}\right)^n,
\eeq
and using \eqn{N424}
then yields \eqn{N417}.

Finally, the case $n=0$ follows from \eqn{N405c} for $n=0$
using \eqn{N422}.
$~~~\bsq$

\paragraph{Remark.} For the case $\pi_{triv}$ on $GL(1)$ Lemma~\ref{le43}
yields 

\begin{eqnarray*}
\tau_n(\pi_{triv}) &= & (-\frac{1}{2})^{n+1} \zeta(n+1, \frac{1}{2}) \\
&=& (-\frac{1}{2})^{n+1} \sum_{m=0}^\infty \frac{1}{(m + \frac{1}{2})^{n+1}} \\
&=&  (-1)^{n+1} \sum_{m=0}^\infty \left( \frac{1}{2m+1} \right)^{n+1}
= (-1)^{n+1} ( 1- \frac{1}{2^{n+1}}) \zeta(n+1).
\end{eqnarray*}
With further work one can deduce
$$
\tau_0(\pi_{triv}) = -\left(\frac{1}{2} \right) 
\left( \log (4 \pi) +  \gamma \right),
$$ 
from which one can deduce  \eqn{N407b}.

%
%
%
\section{Bounds for $S_{\infty}(n, \pi)$}
\hsp
In this and the next section
we consider the order of growth of the terms
$S_\infty(n, \pi)$ and $S_f(n, \pi)$, respectively. 
In this section
we obtain an unconditional result for $S_{\infty}(n, \pi)$.

%
%

\begin{theorem}~\label{th51}
For any irreducible cuspidal (unitary) automorphic representation $\pi$
on $GL(N)$ the quantities $S_{\infty}(n, \pi)$
are real-valued.  There is a constant 
$K(\pi)$ such that for $n \ge K(\pi)$ 
there holds
\beql{501}
S_{\infty}(n, \pi) = \left(\frac{N}{2}\right) n \log n + C_1(\pi)~ n 
+ O \left(N(K(\pi) + 1)\right).
\eeq
Here
\beql{502}
C_1(\pi) =   \frac{N}{2} \left( \gamma - 1 - \log (2\pi) \right) 
+ \frac{1}{2} \log Q(\pi),
\eeq
where $\gamma$ is Euler's constant.
One  can take
\beql{502aa}
K(\pi) = \max_{1\le j \le N}  |\kappa_j(\pi)|^2,
\eeq 
and the  implied constant in the $O$-notation is absolute.
\end{theorem}

Note that $C_(\pi) = C_1(\pi^{\vee})$, since $Q(\pi) = Q(\pi^{\vee})$,
and that  $C_1(\pi)$ does not depend on the archimedean
parameters $\{ \kappa_j(\pi):~ 1 \le j \le N \}$, though these values
appear in the analysis.
For the Li coefficients  we have 
\beql{503aa}
C_1(\pi_{triv}) = 
\frac{1}{2}\left( \gamma  - 1 - \log (2\pi) \right)\simeq -1.1303307. 
\eeq
using $Q(\pi_{triv})=1$.

%
%
%

To begin the proof we introduce the quantities
\beql{506b}
T(n, z) : = \sum_{j=1}^n {\binom{n}{j}}
\left( \frac{1}{2} \right)^j \tilde{\tau}_{j-1}(z).
\eeq
The formula for $S_{\infty}(n, \pi)$
in Lemma ~\ref{le42} can 
be expressed in terms of these for various values of $z$. 
Indeed we obtain, for $n \ge 1$, 
\beql{506a}
S_{\infty}(n, \pi) = \sum_{j=1}^N  T(n, \frac{\kappa_j(\pi^{\vee}) - 1}{2})
+ {\binom{n}{1}} \left( \frac{1}{2} \log Q(\pi) - 
\frac{N}{2} \log \pi \right),
\eeq
using  \eqn{N425}.
We proceed to estimate an individual sum $T(n, z)$. The formulas
for $\tilde{\tau}_j(z)$ give 
$$
T(n, z) = T_1(n, z) + T_2(n,z)
$$ 
in which 
\beql{507}
T_1(n, z) = \frac{n}{2} \frac{\Gamma'}{\Gamma}(1+z)=
{\binom{n}{1}}\frac{1}{2}
( - \gamma + \sum_{m=1}^\infty \frac{z}{m(m+z)}).
\eeq
and 
\beql{N508}
T_2(n, z) := \sum_{j=2}^n {\binom{n}{j}}\frac{1}{2^j} 
\left( (-1)^{j} \sum_{m=1}^{\infty} \frac{1}{(m+z)^j}\right).
\eeq
The second sum converges absolutely and can be rearranged
as
\begin{eqnarray}~\label{N509}
T_2(n, z) &=& \sum_{m=1}^\infty \sum_{j=2}^n {\binom{n}{j}}(-\frac{1}{2})^j
\frac{1}{(m+z)^j} \nonumber \\
&=&  \sum_{m=1}^\infty 
\left( (1 - \frac{1}{2m + 2z})^n - 1 + \frac{n}{2m+2z}\right).
\end{eqnarray}
We divide this sum  into two parts 
$$
T_2(n, z)= T_{20}(n, z) + T_{21}(n, z),
$$
by  splitting
the summation range from $m=[1, n]$ and $m=[n+1, \infty)$
respectively, and  treat these in the following two lemmas.

%
%
%
%
\begin{lemma}~\label{le51}
For any complex number $z$ with  $ \Re(z) \ge -\frac{3}{4}$, the quantity 
\beql{510}
T_{20}(n, z) = \sum_{m=1}^n \left( 
(1 - \frac{1}{2m + 2z})^n - 1 + \frac{n}{2m+2z}\right)
\eeq
satisfies, for all $n \ge |z|^2$, 
\beql{510a}
T_{20}(n, z) = \frac{1}{2}(n\log n) + 
\left(-\frac{1}{2}\frac{\Gamma'}{\Gamma}( 1+ z)
+ \beta_{\infty} + \frac{1}{\sqrt{e}} - 1 \right) n  + O(|z| + 1).
\eeq
in which 
$$
\beta_{\infty} = \int_{1}^{\infty} e^{-\frac{t}{2}} \frac{dt}{t}
\simeq 0.559774.
$$
\end{lemma}
\paragraph{Proof.}
The sum \eqn{510} consists of three subsums, from each
of its terms. The middle subsum in \eqn{510} gives $-n$.
For the third subsum,
\begin{eqnarray}~\label{511}
K(n, z) = \sum_{m=1}^n  \frac{n}{2m+2z} &=& \frac{n}{2} \sum_{m=1}^n 
\left( \frac{1}{m} - \frac{z}{m(m+z)}\right)
\nonumber \\
&=&  \frac{n}{2} \left( \log n + \gamma + O(\frac{1}{n})\right) - 
\frac{n}{2} \left( \frac{\Gamma'}{\Gamma}( 1+ z) + 
\gamma + O(\frac{|z| + 1}{n})\right) 
\nonumber \\
& = & \frac{1}{2} n \log n -  \frac{n}{2}  \frac{\Gamma'}{\Gamma}( 1+ z) +        O(|z| + 1),
\end{eqnarray}
in which we used the partial fraction formula for the digamma 
function and a bound for the error in truncating it at the $n$-th term.

We now  consider the first subsum in \eqn{510}, call it  
\beql{509b}
J(n, z) := \sum_{m=1}^n (1 - \frac{1}{2m + 2z})^n.
\eeq
Setting $z= x+iy$, for real $t\ge 1$ we have 
\beql{511b}
|1 - \frac{1}{2t + 2z}|^2 = 1 + \frac{1 - 4(t+x)}{4( (t+x)^2 +y^2)},
\eeq
The condition  $\Re(z) \ge -\frac{3}{4} $ now implies
that $|1 - \frac{1}{2t + 2z}|^2 \le 1$, so that each term above
is $O(1).$
We  approximate the sum 
by an integral, and  assert that for $n \ge |z|^2$ there holds
\beql{511a}
J(n, z) = \int_1^{n} (1 - \frac{1}{2t + 2z})^n dt + O(|z| + 1).
\eeq
with a constant independent of $z$.
To establish this is we  show that both sum and
integral separately contribute $O(1)$ over  the range 
$1 \le t \le n^{\frac{2}{3}}$,
and then show on the remaining range $n^{\frac{2}{3}} \le t \le n$
that their difference is $O(1)$. 
On the initial range we
use $1- x \le e^{-x}$ for $0 \le x \le 1$ to get from 
\eqn{511b} that, for $t \ge 3$, 
\begin{eqnarray*}
 |1 - \frac{1}{2t + 2z}|^n &\le & \left(1 + \frac{1 - 4(t+x)}{4( (t+x)^2 +y^2)}
\right)^\frac{n}{2} \\
& \le & \exp (\frac{1 - 4(t+x)}{4( (t+x)^2 +y^2)}\frac{n}{2} ) \\
& \le & \exp(-\frac{(t+x)}{ (t+x)^2 +y^2} \frac{n}{4}).
\end{eqnarray*}
Now the function $f(w) = \frac{w}{w^2 + y^2}$ is nonnegative
on the positive real axis, increasing to a maximum at $w=y$ and
decreasing monotonically thereafter. Applied to the
exponential above, we find the the terms for integer $t=m$ with 
$y \le m \le n^{\frac{2}{3}}$ have negative exponent
decreasing in absolute value, so the largest term 
occurs at  the top endpoint,
and is  bounded by $\exp( - c n^{-\frac{1}{3}})$, so the sum
over all these terms is $O(1)$. The first $y$ terms each contribute $O(1)$,
giving $O(|z| + 1)$ in all. This bounds the sum,
and bounding the integral over this range is similar.

On the remaining range $n^{\frac{2}{3}} \le m \le 1$ with
$m \le t \le m+1$ we have  
\begin{eqnarray*}
(1 - \frac{1}{2t + 2z})^n &= & (1 - \frac{1}{2m + 2z})^n
\left( \frac{1 - \frac{1}{2m + 2z}}{1 - \frac{1}{2t + 2z}}\right)^n \\
&= & (1 - \frac{1}{2m + 2z})^n
\left(  1 + O( \frac{1}{|2t+2z|^2} ) \right)^n \\
&= & (1 - \frac{1}{2m + 2z})^n
\left(  1 + O( \frac{n}{|2t+2z|^2} ) \right),
\end{eqnarray*}
where we used $|2t+2z|^2 \ge 2 n.$ Thus,
for $n \ge |z|^2$,  the error between the
sum and the integral is bounded by an absolute constant times 
\begin{eqnarray}
\int_{n^{2/3}}^{n} (1- \frac{1}{2t + 2z})^n \frac{n}{|2t+2z|^2} dt
& \le & O( 1).
\end{eqnarray}
This last estimate is obtained by observing that the integrand is 
maximized (in absolute value) at its top endpoint, where it 
is $O( \frac{1}{n})$. We conclude that \eqn{511a} holds.

We now integrate by parts on the right side of \eqn{511a} to
obtain 
$$J(n, z)  = J_1(n, z) - J_2(n,z) + O (|z|+1),
$$
 in which
\beql{514a}
J_1(n, z) :=  t(1- \frac{1}{2t+2z})^n \vert_{t=1}^n = n(1- \frac{1}{2n+2z})^n
- (1 -\frac{1}{2+2z})^n,
\eeq
and
\beql{514b}
J_2(n, z) = n \int_{1}^n \frac{2t}{(2t+2z)^2} (1- \frac{1}{2t+2z})^{n-1} dt.
\eeq
The first term on the right side is estimated for $n > |z|^2$ by
\beql{N522}
J_1(n,z)  = n \exp( -\frac{n}{2n + 2z})
\left(1 +   O(\frac{|z| + 1}{n}\right) +O(1)
 = \frac{1}{\sqrt{e}} n  + O ( |z| + 1),
\eeq
For the  remaining integral $J_2(n, z)$, we show that
$\tilde{J}_2(n, z) = \frac{1}{n}J_2(n, z)$ 
approaches a limit as $n \to \infty$, using
the fact that $(1 - \frac{1}{u})^u$ approaches $ e$ as
$\Re(u) \to \infty$.
We rescale it with the variable change $u=nt$ to obtain  
\begin{eqnarray*}
\tilde{J}_2(n, z) & = & \int_{\frac{1}{n}}^1 \frac{2nu}{(2nu+z)^2}
\left( 1 - \frac{1}{2(nu+z)}\right)^{n-1} n du \\
& = & \frac{1}{2} \int_{\frac{1}{2n}}^1 \frac{u}{(u + \frac{2z}{n})^2}
\left[ \left(  1 - \frac{1}{2(nu+z)}\right)^{2(nu+z)}\right]^{\frac{n-1}{2(nu+z)}}
du \\
&=& \frac{1}{2}\int_{\frac{1}{n}}^1 \frac{1}{u} e^{-\frac{1}{2u}} du + 
O \left(\frac{|z|+1}{n} \right) \\
&=& \frac{1}{2}\int_{0}^1 \frac{1}{u} e^{-\frac{1}{2u}} du + 
O \left(\frac{|z|+1}{n}\right)
\end{eqnarray*}
Now we set 
$$
\beta_{\infty} := \int_{0}^1 e^{-\frac{1}{2u}} \frac{du}{u} = 
\int_{1}^{\infty} e^{- \frac{t}{2}} \frac{dt}{t},
$$
and we have obtained
$$
J_2(n,z) = \frac{1}{2}\beta_{\infty} n + O ( |z| + 1).
$$
Combining all these estimates gives, for $n \ge |z|^2$, 
\begin{eqnarray*}
T_{20}(n, z) &= & J_1(n, z) - J_2(n, z) -n + K(n,z) + O(|z|+1)\\
&=& ( \sqrt{e}-1)n +
\frac{1}{2} n \log n -  \frac{n}{2}\frac{\Gamma'}{\Gamma}(1+z) 
- \frac{1}{2}\beta_{\infty} n +   
O(|z| + 1)\\
& = & \frac{1}{2} n \log n + \left( -\frac{1}{2}\frac{\Gamma'}{\Gamma}(1+z)-
\frac{1}{2}\beta_{\infty} + \frac{1}{\sqrt{e}} - 1 \right) n
+ O (|z| + 1),
\end{eqnarray*}
the desired estimate.
$~~~\bsq$ \\

%
%
%
%
\begin{lemma}~\label{le52}
For $\Re(z) \ge -1$ and $n \ge |z|+2$, the quantity 
$$
T_{21}(n, z) = \sum_{m=n+1}^{\infty} \left( 
(1 - \frac{1}{2m + 2z})^n - 1 + \frac{n}{2m+2z}\right)
$$
satisfies
\beql{523}
T_{21}(n, z) = (\frac{1}{2} - \frac{1}{\sqrt{e}} 
+\frac{1}{2}\alpha_{\infty}) n + O(|z|+1),
\eeq
in which
\beql{524}
\alpha_{\infty} = 
\int_{0}^{1} (1 - e^{-\frac{1}{2}t})\frac{dt}{t} \simeq 0.443842.
\eeq
The implied  $O$-constant is absolute.
\end{lemma}

\paragraph{Proof.} 
In what follows we  assume $n \ge |z| + 2$.
We approximate the sum by an integral
\beql{525}
T_{21}(n, z) = \int_{n}^{\infty} \left( 
(1 - \frac{1}{2t + 2z})^n -1 + \frac{n}{2t + 2z} \right) dt + O(1).
\eeq
To justify  the error term, 
we observe that for $n \le m \le t \le m+1$ there holds
$$
(1 - \frac{1}{2t + 2z})^n =  \left(1 - \frac{1}{2m + 2z}+ 
O(\frac{1}{|(2m+2z)(2t+ 2z)|}) \right)^n = (1 - \frac{1}{2m + 2z})^n
\left(1 + O (\frac{n}{(t+|z|)^2}\right)
$$
where we used $\Re(z) \ge -1$ and $n \ge 2$.
We also use 
$$
|\frac{n}{t + z} - \frac{n}{m+z}| \le |\frac{n(m-t)}{|(t+z)(m+z)}|
\le \frac{n}{(m + |z| - 1)^2},
$$
Viewing the sum as an integral of a step function, we have
bounded the difference between the integrands at $t$ by an
absolute constant times $\frac{n}{(t+|z|-1)^2}$, and 
for $n \ge 3$ we have
$$
\int_{n}^{\infty} \frac{n}{(t+|z|-1)^2} dt \le 
n \int_{n-1}^{\infty} \frac{1}{t^2} dt = O(1).
$$
and \eqn{525} follows.

Making the change of variable $u= \frac{1}{2t+2z}$ we obtain
\beql{526}
T_{21}(n, z) = \frac{1}{2}\int_{0}^{\frac{1}{2n+2z}} 
[ (1 - u)^n - 1 + nu] \frac{du}{u^2} + O(1).
\eeq
This integral is a contour integral 
(since $z$ is complex) but the answer is independent
of the contour, since the integrand is analytic.
 An integration by parts yields
$T_{21}(n, z) = K_1(n, z) - K_2(n, z) + O(1)$
in which
\begin{eqnarray}~\label{527a}
K_1(n, z) &=&  -\frac{1}{2u} [(1 - u)^n - 1 + nu] |_{u=0}^{\frac{1}{2n+2z}}
\nonumber \\
&=& -(n+z)\left((1- \frac{1}{2n+2z})^n - 1 + \frac{n}{2n+2z}\right),
\end{eqnarray}
and 
\beql{527b}
K_2(n,z) := \frac{1}{2}\int_{0}^{\frac{1}{2n+2z}}
[ -n(1 - u)^{n-1} + n] \frac{du}{-u} = 
\frac{n}{2}\int_{0}^{\frac{1}{2n+2z}}[ (1- u)^{n-1} - 1]\frac{du}{u}.
\eeq 
To estimate $K_1(n,z)$ we write 
$z = x+iy$, with $x \ge -1$, to obtain  
$$
1- \frac{1}{2n+2z}= ( 1- \frac{2n+ 2x}{4(n+x)^2 + 4y^2})
+ i \frac{2y}{ 4(n+x)^2 + 4y^2} = r e^{i \theta},
$$
in which $r = 1 - \frac{1}{2n} + O (\frac{1}{n^2})$ and 
$\theta = O( \frac{|z| + 1}{n^2}).$ It follows that
$r^n e^{i n \theta} = e^{-\frac{1}{2}} + O (\frac{|z| + 1}{n})$,
and we obtain 
\begin{eqnarray}~\label{527c}
K_1(n, z) &=& -n[(1- \frac{1}{2n+2z})^n  - \frac{1}{2}] +O(|z|+1) 
\nonumber \\
&=& 
-n [r^n e^{i n \theta} -  \frac{1}{2}] + O(|z|+1)
\nonumber \\
& = & n\left( \frac{1}{2}- \frac{1}{\sqrt{e}} + O (\frac{|z|+1}{n})\right) 
+ O(|z|+1).
\end{eqnarray}
Expanding the integral $K_2(n, z)$ in powers of $u$ yields
\begin{eqnarray}~\label{528}
K_2(n, z) &=& \frac{n}{2} \left(\sum_{j=1}^n \int_{0}^{\frac{1}{2n+2z}} 
(-1)^j {\binom{n}{j}} u^{j-1} du\right) 
\nonumber \\
&=& \frac{n}{2} \left(
\sum_{j=1}^n \frac{(-1)^j}{j} {\binom{n}{j}}
(\frac{1}{2n + 2z})^j\right).
\end{eqnarray}
We set 
$$
\alpha_n(z) = \sum_{j=1}^n \frac{(-1)^j}{j} 
{\binom{n}{j}}(\frac{1}{2n + 2z})^j,
$$
and  
\beql{528a}
\alpha_{\infty} := \sum_{j=1}^\infty \frac{(-1)^{j+1}}{j} \frac{1}{j! 2^j}
= \int_{0}^{\frac{1}{2}} (1 - e^{-v}) \frac{dv}{v}
= \int_{0}^{1} (1 - e^{-\frac{1}{2}t}) \frac{dt}{t}.
\eeq
We assert that  for $\Re(z)\ge -1$, $n \ge 2$, 
\beql{529}
|\alpha_n(z) + \alpha_{\infty}| = O(\frac{|z|+1}{n}).
\eeq
where the constant in the O-notation is independent of $z$.
Indeed we note that
$$
{\binom{n}{j}}(\frac{1}{2n + 2z})^j 
=  \frac{1}{j! 2^j} (\frac{n}{n+z})(\frac{n-1}{n+z}) \cdots
(\frac{n-j+1}{n+z}),
$$
which since $\Re(z) \ge -1$ yields 
$$
|{\binom{n}{j}}(\frac{1}{2n + 2z})^j| \le  \frac{2}{j! 2^j}.
$$
and also for  $0 \le j \le 2 \log n$ that
$$
|{\binom{n}{j}}(\frac{1}{2n + 2z})^j- \frac{1}{j! 2^j}|=
O \left(\frac{|z|+1}{n} \right).
$$
Combining these estimates with 
$$
|\alpha_{\infty} + \sum_{j =1}^{2 \log n} \frac{ (-1)^j}{j} \frac{1}{j ! 2^j}|
\le O(\frac{1}{n}),
$$
yields \eqn{529}, so that
$K_2(n, z) = -\frac{1}{2} \alpha_{\infty} n + O (|z| +1)$.

Combining all the estimates above gives
\begin{eqnarray}~\label{530}
T_{21}(n, z) &=& K_1(n, z) - K_2(n,z) + O(1))
\nonumber \\
& = & (\frac{1}{2} - \frac{1}{\sqrt{e}} +\frac{1}{2}\alpha_{\infty}) n + O (|z| + 1), \nonumber
\end{eqnarray}
the desired estimate.
$~~~\bsq$

%
%

\paragraph{Proof of Theorem~\ref{th51}.}
Lemma~\ref{le42} gives 
\begin{eqnarray*}\label{531}
S_{\infty}(n, \pi) & =& \sum_{j = 1}^n {\binom{n}{j}}\tau_{j-1}(\pi)
\nonumber \\
&=& \left(\sum_{k=1}^N \sum_{j = 1}^n {\binom{n}{j}} (-1)^j 
\tilde{\tau}_{j-1}(\frac{\kappa_j(\pi)-1}{2})(\frac{1}{2})^j\right)
+ {\binom{n}{1}}\left( \frac{1}{2} \log Q(\pi) -  \frac{N}{2} \log \pi \right)
\nonumber \\
&=& 
\left(\sum_{k=1}^N T\left(n, \frac{\kappa_k(\pi) -1}{2}\right) \right) 
+\left(\frac{1}{2} \log Q(\pi) - \frac{N}{2} \log \pi\right) n
\nonumber \\
&=& \left( \frac{1}{2}\log Q(\pi) - \frac{N}{2} \log \pi +
 \frac{1}{2}\sum_{j=1}^N 
\frac{\Gamma'}{\Gamma} \left(\frac{\kappa_k(\pi) +1}{2}\right) \right) n  +
\sum_{k=1}^N T_2\left(n, \frac{\kappa_k(\pi) -1}{2}\right).
\end{eqnarray*}
where we made use of the decomposition 
$T(n, z) = T_1(n, z) + T_2(n, z)$ and the formula \eqn{507}
for $T_1(n, z)$.

Now we suppose that 
$n \ge K(\pi) :=\max \{ |\kappa_k(\pi)|^2: 1 \le k \le N\}$ 
and note by Theorem~\ref{th21}(2)
that each $\Re( \frac{\kappa_k(\pi) +1}{2}) > -\frac{3}{4}$,  
so that we can
apply Lemma~\ref{le51} and \ref{le52}. We  obtain
\begin{eqnarray}
T_2(n,~ \frac{\kappa_k(\pi) -1}{2}) &= &
\frac{1}{2} n \log n + 
\left(- \frac{1}{2} \frac{\Gamma'}{\Gamma}(\frac{\kappa_k(\pi) + 1}{2}) +
\frac{1}{\sqrt{e}} - 1 \right) n \nonumber \\
&&~~~~~~~~ +~~
(\frac{1}{2} - \frac{1}{\sqrt{e}} + 
\frac{1}{2}\alpha_{\infty}- \frac{1}{2}\beta_{\infty})n + O(K(\pi) + 1))
\nonumber \\
&=& 
\frac{1}{2} n \log n + \left( 
-\frac{1}{2}\frac{\Gamma'}{\Gamma}(\frac{\kappa_k(\pi) + 1}{2}) + 
\frac{1}{2}(\alpha_{\infty} -\beta_{\infty}- 1)\right) n + O(K(\pi) +1). \nonumber
\end{eqnarray}
Substituting this in the previous formula  yields 
\begin{eqnarray}~\label{532}
S_{\infty}(n, \pi) &= & \left(\frac{N}{2}\right)  n \log n 
+ \left(  \frac{N}{2}(\alpha_{\infty}-\beta_{\infty}- 1 -  \log \pi) +
\frac{1}{2} \log Q(\pi) \right) ~n
\nonumber \\
&& ~~~~~~~~+~~ O \left( N (K(\pi) + 1)\right). 
\end{eqnarray}
To complete the proof we must establish 
the identity
\footnote{This identity was found by comparison of 
the author's original formula \eqn{532} with 
an  asymptotic formula for the Li 
coefficients (under RH) given by 
Voros~\cite[eqn. (11)]{Vo04}.}
\beql{503ab}
\alpha_{\infty} - \beta_{\infty}= \gamma - \log 2,
\eeq
where $\gamma$ is Euler's constant.
We start from 
$$
\alpha_{\infty} - \beta_{\infty}= 
\int_0^{\frac{1}{2}} (1 - e^{-t}) \frac{dt}{t}
-\int_{\frac{1}{2}}^{\infty} e^{-t} \frac{dt}{t},
$$
The assertion  \eqn{503ab} is  a special case
of an  identity valid for each  $w >0$, that
\beql{534}
\int_0^{w} (1 - e^{-t}) \frac{dt}{t} - \int_{w}^\infty e^{-t}\frac{dt}{t}
= \gamma + \log w,
\eeq
on taking $w = \frac{1}{2}.$ The
reviewer observes the following proof of \eqn{534}.
Denoting the left side of this
equation by $S(w)$, an integration by parts  (taking $dv= \frac{dt}{t}$)
yields
\beql{534a}
S(w)= \log w - \int_{0}^{\infty} e^{-t}\log t dt = \log w - \Gamma^{'}(1),
\eeq
and the result follows using $\Gamma^{'}(1)= - \gamma$.
This method was followed by Barnes \cite{Ba01} in establishing 
the identity
\beql{535}
\gamma = \int_0^1 (1 - e^{-t} - e^{-\frac{1}{t}}) \frac{dt}{t}.
\eeq
$~~~\bsq$.

\paragraph{Remark.} Computational evidence of Maslanka
(private communication) indicates that
$S_{\infty}(n, \pi)$ has a full  asymptotic expansion in
inverse powers $n^{-k}$ with $k \ge 0$. The further coefficients
of such expansions will likely depend on the 
archimedean parameters 
$\{\kappa_j(\pi): 1 \le j \le N\}$,
unlike the first two coefficients appearing in Theorem~\ref{th51}.
%
%
%

\section{Bounds for $S_{f}(n, \pi)$}
\hsp
We obtain a bound for $S_{f}(n , \pi)$
in terms of the 
{\em incomplete  Li coefficient to height $T$,} defined by
\beql{N601}
\lambda_n(T, \pi) := \sum_{{{\rho \in Z(\pi)}\atop{|\Im(\rho)| \le  T}}}
1 - \left( 1 - \frac{1}{\rho}\right)^n,
\eeq 
where $T$ is a cutoff parameter.

\begin{theorem}~\label{th61}
For any irreducible cuspidal (unitary) automorphic representation
on $GL(N)$  there holds
\beql{N602}
S_{f}(n, \pi) = \lambda_n(  \sqrt{n}, \pi^{\vee})  + 
O\left(  \sqrt{n} \log n \right),
\eeq
in which the implied constant in the $O$-notation depends
on  $\pi$. If the Riemann hypothesis holds
for $L(s, \pi)$ then
\beql{N602aa}
\lambda_n(\sqrt{n}, \pi^{\vee}) = O \left( \sqrt{n} \log n \right).
\eeq
\end{theorem}

\paragraph{Proof.}
We use a contour integral argument. We introduce the kernel function
\beql{N602b}
k_n(s) := (1 + \frac{1}{s})^n  - 1 = \sum_{j=1}^n 
{\binom{n}{j}}\left(\frac{1}{s}\right)^j,
\eeq
If $C_1$ is a contour enclosing the point $s=0$ counterclockwise
on a circle of small enough positive radius $R$, the residue theorem gives
\beql{N602a}
\frac{1}{2 \pi i} \int_{C_1} k_n(s)
\left( - \frac{L'}{L}(s+1, \pi) \right) ds = 
\sum_{j=1}^n  {\binom{n}{j}} \eta_{j-1} = S_f(n, \pi).
\eeq
The residue comes entirely
 from the singularity at $s=0$, as no other
singularities lie inside the contour. 

We deform the contour to the counterclockwise oriented rectangular
contour $C_2(n)$ consisting of 
 vertical lines with real part  $\Re(s) = \sigma_0$
and $\Re(s) = \sigma_1$
where we will choose $-3< \sigma_0 < -2$,
and 
$\sigma_1 = 2\sqrt{n}$,
and horizontal lines at 
$\Im(s) = \pm T$,
where we will choose $T=\sqrt{ n} + \epsilon_n$,
for some $0 < \epsilon_n < 1$. 
The integrand has simple poles
at the zeros of $L(s, \pi)$,
and some of them will now lie inside the contour. 
The residue theorem gives
$$
I_2(n) := \int_{C_2(n)}  k_n(s)
\left( - \frac{L'}{L}(s+1, \pi) \right) ds =
S_{f}(n, \pi) + \sum_{{\rho \in Z(\pi)}\atop{ |\Im(\rho)| < T }}
\left(1 + \frac{1}{\rho - 1}\right)^n - 1
+ I_2^{{triv}}(n),
$$
in which $I_2^{{triv}}(n)$ evaluates the residues coming from the
trivial zeros of $L(s, \pi)$ which will satisfy 
$ -3 \le \Re(s) < \frac{1}{2}$. The trivial zeros are associated
to the archimedean factors and their location is dictated by
the values $\kappa_j(\pi)$, which 
by Theorem~\ref{th21} satisfy $\Re(\kappa_j(\pi)) > -\frac{1}{2},$
so that all trivial zeros satisfy $\Re(\rho) < \frac{1}{2}.$ 
We have 
\beql{N603}
| 1+ \frac{1}{z -1}| = |\frac{z}{z-1}|< 1 ~~~\mbox{when}~~~ 
\Re(z) < \frac{1}{2},
\eeq
and using this we conclude that the residue at each trivial
zero is $O(1)$. Since the  trivial zeros fall in $N$ arithmetic
progressions with spacing $2$, there 
are $O(N)$ such zeros in the interval
$ -3 \le \Re(s) < \frac{1}{2}$ we conclude that the
trivial zero contribution is  
$I_2(n, triv) = O(1)$, with $O$-constant depending on $\pi$.

Using the symmetry $\rho \mapsto 1 - \bar{\rho}$ of $Z(\pi)$,
which does not change $|\Im(\rho)|$, we can rewrite the sum over
$Z(\pi)$ in terms of 
$$
\left( \frac{1-\bar{\rho}}{- \bar{\rho}}\right)^n  -1
= \left( \frac{\bar{\rho} - 1}{\bar{\rho}}\right)^n -1.
$$
We observe that $\bar{\rho} \in Z(\pi^{\vee})$ so that \eqn{N602a}
can be rewritten
\beql{N604}
I_2(n) = S_{f}(n, \pi) - \lambda_n(T, \pi^{\vee}) + O(1).
\eeq
The main part of the theorem will be to establish
$$
I_2(n) = O \left( \sqrt{n} \right).
$$
The theorem will follow from the two assertions that
$I_2(n) = O \left( \sqrt{n}\right)$
and 
$$
|\lambda_n( \sqrt{n}, \pi^{\vee}) - 
\lambda_n(T, \pi^{\vee})| = O (\log n).
$$
The second assertion  follows from the 
observation that $|T - \sqrt{n}| <1$,
that there are $O(\log n)$ zeros in an interval of
length one at this height (by \eqn{N208}), and that for each zero
$\rho= \beta + i \gamma$ with
$\sqrt{n} \le |\Im(\rho)| < \sqrt{n} +1$ 
there holds
\begin{eqnarray*}
|\left( \frac{\rho-1}{\rho} \right)^n| 
&= & |1-\frac{\beta-i\gamma}{\beta^{2}+ \gamma^{2}}|^n \\
& \le & |\left(1- \frac{\beta}{\beta^{2} + \gamma^{2}}\right)^{2}
+ \left(\frac{\gamma}{\beta^{2}+ \gamma^{2}}\right)^{2}|^{n/2}\\
&\le & |1 + \frac{1}{n}|^{n/2} \le 2.
\end{eqnarray*}

It remains to bound $I_2(n)$. We will use the following
estimate for $\frac{L'}{L}(s, \pi)$ in a region
including the  critical strip.

%
%

\begin{lemma}~\label{le61}
For $-4 \le \Re(s) < 4$ and $s= \sigma + it$ there holds
\beql{N607}
 \frac{L'}{L}(s, \pi) = \sum_{\{\rho: |\Im(\rho - s)| < 1\}} 
\frac{1}{s- \rho} + O \left(\log \fq(\pi,s) \frac{L'}{L}(2, \pi)\right),
\eeq
in which $\fq(s, \pi)= Q(\pi) \prod_{i=1}^N( |s + \kappa_j(\pi)| + 3)$ 
is the analytic conductor of $\pi$, and 
the constant in the $O$-notation is absolute.
\end{lemma}

\paragraph{Proof.} 
The classical proof for Dirichlet $L$-functions  
\cite[Chap. 15, p. 102]{DM00} generalizes to give \eqn{N607}, following
Chapter 5 of \cite{IK04}.
If we know the Ramanujan conjecture holds for $L(s, \pi)$,
then we can control the size of $\frac{L'}{L}(2, \pi)$ in the 
remainder term.
$~~~\bsq$. \\

We now choose the parameters $\sigma_0$ and $T$ appropriately
to avoid poles of the integrand.
Since the trivial zeros of $L(s, \pi)$ fall in $N$ arithmetic
progressions with spacing $2$, we may choose $\sigma_0$ so that the contour
avoids coming within $\frac{1}{2N}$ of
any trivial zero. Similarly we can choose 
$T = \sqrt{n} + \epsilon_n$ with $0 \le \epsilon_n \le 1$
so that the horizontal lines do not approach closer than
$O(\log n)$ to any zero of $L(s, \pi)$. It then follows
from Lemma~\ref{le61} that on the horizontal line in the 
interval $-3 \le \Re(s) \le 4$ we have
\beql{N609}
 |\frac{L'}{L}(s+1, \pi)| = O( (\log T)^2).
\eeq

The Euler product for $L(s, \pi)$
converges absolutely for $\Re(s) > 1$ (by Theorem~\ref{th21}(1))
hence the Dirichlet
series for $\frac{L'}{L}(s, \pi)$ converges
absolutely for $\Re(s) > 1.$ More precisely if 
$\frac{L'}{L}(s, \pi)= \sum_{m=2}^{\infty} b(m, \pi) m^{-s}$
then for $\sigma >1$,
$$|\frac{L'}{L}|(\sigma) = \sum_{m=2}^{\infty} |b(m, \pi)| m^{-\sigma} 
< \infty.
$$ 
Since the Dirichlet series
for $\frac{L'}{L}(s, \pi)$ has no constant term, we obtain 
for $\Re(s) = \sigma > 2$ the bound 
\beql{N610}
|\frac{L'}{L}(s, \pi )| \le |\frac{L'}{L}|(\sigma, \pi)| \le       
2^{-(\sigma - 2)} |\frac{L'}{L}(2, \pi)|.
\eeq

We now consider the integral $I_2(n)$ on the vertical segment (I)
having  $\Re(s) = \sigma_1= 2 \sqrt{n}$, call it $I_2^{(I)}(n)$. We have
$$
|(1 -\frac{1}{s})^n - 1| \le (1+ \frac{1}{\sigma_1})^n +1
\le (1 + \frac{1}{2\sqrt{n}})^n \le \exp(\frac{1}{2} \sqrt{n}) < 2^{\sqrt{n}}.
$$
Now the estimate \eqn{N610} gives on (I) that
$$
|\frac{L'}{L}(s, \pi )| \le C_0 2^{-2(\sqrt{n} +2)}.
$$
so that the integrand is bounded above by $C_1 2^{- \sqrt{n}}$.
The length of the contour is $O(\frac{n}{\log n})$ so
we obtain the estimate $|I_2^{(I)}(n)| = o(1)$.

We next consider the integral $I_2(n)$ on the two horizontal
segments, call them (II) and (IV). It  suffices to bound
$I_2^{(II)}(n)$, the treatment for (IV) being identical.
Let $s = \sigma +it$ be a point on (II).
We have  $T \ge \sqrt{ n}$ so that
\begin{eqnarray}~\label{N617}
| 1 + \frac{1}{s}| &= & | (1 + \frac{\sigma}{\sigma^2 + T^2}) + 
\frac{iT}{\sigma^2 + T^2})| \nonumber \\ 
&\le &  \left( (1 + \frac{\sigma}{\sigma^2 + T^2})^2 
+ \frac{1}{\sigma^2 + T^2}  \right)^{\frac{1}{2}}  \nonumber \\ 
& \le & 1 + \frac{\sigma+1}{\sigma^2 + T^2}.
\end{eqnarray}
By hypothesis $T^2 \ge  n$ so for $-3 \le \sigma \le 3$
we have 
$$
|k_n(s)| \le ( 1 + \frac{4}{16 + n})^n +1 \le
e^{4} + 1=  O \left( 1 \right)
$$
On this interval Lemma~\ref{le61} gives
$$
|\frac{L'}{L} (s, \pi)| = O ((\log T)^2) = O ((\log n)^2)
$$
since we have chosen the ordinate $T$ to stay away from
zeros of $L(s, \pi)$. 
Now we step across the interval (II) towards
the right, in segments of
length $1$, starting from $\sigma= 3$. At the initial
point we have 
$|\frac{L'}{L} (s, \pi)|= O(1)$ since we are in the
absolute convergence region, and the estimate
\eqn{N610} gains a factor of $2$ going one unit to the
right. In comparison 
\begin{eqnarray}~\label{N620}
\left|\frac{k_n(s+1)+1}{k_n(s)+1}\right| &= & 
| \frac{1 + \frac{1}{\sigma +1 + iT}}{1 + \frac{1}{\sigma + iT}}|^n
 \nonumber \\
&\le & |1 + (\frac{1}{\sigma +1 + iT}-\frac{1}{\sigma + iT})|^n 
 \nonumber \\
&\le & \left(1 + \frac{1}{|\sigma + 1 + iT|~|\sigma + iT|} \right)^n 
\nonumber \\
&\le & \left(1 + \frac{1}{T^2}\right)^n \le e.
\end{eqnarray}
We obtain an upper bound for $|k_n(s) \frac{L'}{L}(s, \pi)|$
that decreases geometrically at each step, and after $O(\log n)$
steps it becomes $O(1)$ so we obtain the upper bound
$$
|I_2^{(II)}(n)| = O( \sqrt{n} + (\log n)^3) = O (\sqrt{n}).
$$
A similar bound holds for $|I_2^{(IV)}(n)|.$

For the remaining vertical segment (III) with $\Re(s) = \sigma_0$,
we have that the  kernel function $|k_n(s)|= O(1)$ on the line
segment (III) by \eqn{N603},
and  $|\frac{L'}{L}(s, \pi)| = O (\log (|s| + 2))$
using Lemma~\ref{le61}. Since the segment (III)
has length $O(\sqrt{n})$ we obtain the bound
\beql{N615a}
 |I_2^{(III)}(n)| = O(\sqrt{n} \log n).
\eeq
Using this bound  suffices to establish  \eqn{N602},
with the remainder term $O( \sqrt{n} \log n )$.

Totalling all these bounds above gives
$$
S_{f}(n, \pi) = \lambda_n(T, \pi^{\vee}) + O \left( \sqrt{n} \log n \right),
$$
with $T = \sqrt{ n} + \epsilon_n$.

Now suppose  that the Riemann hypothesis holds for $\xi(s, \pi)$.
Then we have 
$$
|1 - \frac{1}{\rho -1}| = |\frac{\rho}{\rho -1}| =
|\frac{-\frac{1}{2} + i \gamma}{\frac{1}{2} + i \gamma} |=1
$$
It follows that each zero contributes a term of absolute value
at most $2$ to the incomplete Li coefficient $\lambda_n(T, \pi)$,
and we obtain 
\beql{666}
\lambda_n(T, \pi) = O(T \log T+ 1),
\eeq
on using the zero density estimate in Theorem~\ref{th21}(4).
Here the constant in the $O$-notation depends on $\pi$.
In particular
$$
\lambda_n(\sqrt{ n}, \pi^{\vee}) = \overline{\lambda_n(\sqrt{n}, \pi)}
= O (\sqrt{n} \log n),
$$
as required.
$~~~\bsq$

%
%
%
\section{Interpolation Function for  Li Coefficients} 
\hsp
We construct an entire 
function that interpolates the Li coefficients at integer
values.

\begin{theorem}~\label{th71}
(1) For any irreducible cuspidal unitary automorphic representation
$\pi$ over $GL(N)$ over $\QQ$, there exists an 
entire function $F_{\pi}(z)$
of order one and exponential type having the 
two properties:

(i). It interpolates
the generalized Li coefficients $\lambda_n(\pi)$
at integer values, i.e.
\beql{N701}
F_{\pi}(n) = \lambda_n(\pi),~~~ n \in \ZZ.
\eeq

(ii) It is real-valued on the imaginary axis, with 
\beql{N702}
F_{\pi}(-\bar{z}) = \overline{F_{\pi}(z)}.
\eeq

\noindent (2) If the Riemann hypothesis holds
for $L(s, \pi)$, 
then there exists a unique
function  $F_{\pi}(z)$ of
exponential type at most $\pi$ that satisfies  (i), and is such that
$$
F_{\pi}^{\ast}(z) : = F_{\pi}(z) - e(\frac{1}{2},\pi)(1-\cos \pi z)
$$
is of exponential type strictly less than $\pi$. Here $e(\frac{1}{2},\pi)$
is the order of the zero of $\xi(s,\pi)$ at $s=\frac{1}{2}$.
On the real axis, $F_{\pi}(z)$ satisfies the bound
\beql{N702a}
|F_{\pi}(x)| \le C (|x|+2) \log (|x| +2), 
\eeq
for some constant $C$ depending on $\pi$.
\end{theorem}

\paragraph{Proof.}
In this proof we regard $\pi$ as fixed, with zeros 
$\rho$ drawn from  the multiset $Z(\pi)$.

(1). We define
\beql{702}
(1 - \frac{1}{\rho})^z := \exp \left( z \log( 1 - \frac{1}{\rho}) \right),
\eeq
in which we take the principal branch of the logarithm,
with argument $-\frac{\pi}{2} < \Im(\log z) \le \frac{\pi}{2}.$

We also have the binomial expansion
\beql{703}
(1 - \frac{1}{\rho})^z = \sum_{j=0}^{\infty} 
(-1)^j{\binom{z}{j}} \frac{1}{\rho^j}.
\eeq
When $|\rho| > 1$ this expansion converges absolutely
for all $z \in \CC$, since the terms in the series
eventually decay geometrically in absolute value.

We now let $R >0$ and define
\beql{704}
F_R(z) = \sum_{|\rho| \le   R} \left( 1 - (1 - \frac{1}{\rho})^z\right)  +
  (-z) \left(\sum_{\{|\rho| > R}{}^{'} \frac{1}{\rho}\right) +
\sum_{j=2}^{\infty} (-1)^j{\binom{z}{j}} \left(
\sum_{|\rho| > R} \frac{1}{\rho^j} \right),
\eeq
Here  
\beql{705}
\sum_{|\rho| > R} {}^{'} \frac{1}{\rho} :=
\lim_{T \to \infty} \sum_{R < |\rho| \le T} \frac{1}{\rho},
\eeq
and this sum is $*$-convergent by Lemma~\ref{Nle21}.
For notational convenience we define the  {\em partial power sums}
\beql{705a}
 \sigma_j(R, \pi) := \sum_{|\rho| > R} \frac{1}{\rho^j}.
\eeq
These sums converge absolutely for $j \ge 2$ and are
$*$-convergent for $j=1$, by Lemma~\ref{Nle21}.
We can rewrite
\beql{706}
F_R(z) = \sum_{|\rho| \le   R} (1 - \frac{1}{\rho})^z
 - \sigma_1(R, \pi) z
+ \sum_{j=2}^{\infty} (-1)^j{\binom{z}{j}} \sigma_j(R, \pi)
\eeq
The zero-counting estimate in Theorem~\ref{th21}(4) implies
there are at most $O(\log T)$ zeros in an interval of length
$1$ at height $T$, where the $O$-constant depends on the
representation $\pi$. Therefore we obtain, for $j \ge 2$
\begin{eqnarray}~\label{707}
|\sigma_j(R, \pi)| & \le &  \sum_{m=R}^{\infty} 
\sharp  \{ m \le |\rho| < m+1 \} \frac{1}{m^j} 
\nonumber \\
 & \le & \sum_{m=R}^{\infty} C \frac{\log m}{m^j} 
= O \left( (\log R)^2  \frac{1}{R^{j-1}} \right).
\end{eqnarray}
where the $O$-constant depends on $\pi$.
For the term $j=1$ we similarly obtain a $*$-convergence estimate
\beql{708}
|\sigma_1(R, \pi)| = O (\frac{ (\log R)^2}{R})
\eeq
by partial summation, cancelling the zeros at height $T$
against those at height $-T$. 
We therefore have, on the disk $|z| \le \sqrt{R}$, that

\begin{eqnarray}~\label{709}
| {\binom{z}{j}} \sigma_j(R, \pi)| &\le & |{\binom{z}{j}}| |\sigma_j(R, \pi)|
\nonumber \\
&\le & C_1 \left( \prod_{k=1}^j \frac{ \sqrt{R} + k}{k}\right) (\log R)^2 
R^{-j +1} \nonumber \\
&\le & C_1 \left(  \prod_{k=1}^j \frac{1+ \frac{k}{\sqrt{R}}}{k} \right) 
(\log R)^2  R^{-\frac{j}{2} +1}  \nonumber \\
&\le & C_2 (\log R)^2  R^{-\frac{j}{2} +1}.
\end{eqnarray}
On summing over $j \ge 1$ it follows  that the terms on 
the right side of  \eqn{706} converges uniformly
on the entire disk $|z| < \sqrt{R}$ to an analytic function
$F_R(z)$. 

We also obtain from this a 
bound on the maximum modulus of  $F_R(z)$ 
on $|z| \le R$. Set 
\beql{710}
\beta(\pi) := \sup_{\rho \in Z(\pi)} |\log(1- \frac{1}{\rho})|,
\eeq
and observe $\beta(\pi)$  is finite because $1 \notin Z(\pi)$,
and $|1 -\frac{1}{\rho})| \to 0$ as $|\rho| \to \infty.$ 
then we obtain, for $|z| \le \sqrt{R} \ge 2$, that
\begin{eqnarray}~\label{711}
|F_R(z)| &\le & \sum_{|\rho| \le R} e^{\sqrt{R} \log(1- \frac{1}{\rho})}
+ O \left( (\log R)^2 \frac{1}{1- R^{-\frac{1}{2}}} \right)
\nonumber \\
&=& C_3 \left(R\log R\right) e^{\beta(\pi)\sqrt{R}}
\end{eqnarray}
in which the $O$-constant depends on $\pi$.

We assert that all  $F_R(z)$ represent the same analytic
function. If $R_1 < R_2$ and we consider the functions on
the domain $|z| < \sqrt{R_1}$ we can expand the individual
terms
$(1 - \frac{1}{\rho})^z$ with $R_1 < |\rho| \le R_2$ 
in $F_{R_2}(z)$ and combine them term-by-term with
the binomial expression, using
$$
\sigma_j(R_1, \pi) = \sigma_j(R_2, \pi) + 
\sum_{R_1 < |\rho| \le R_2} \frac{1}{\rho_j},
$$
Since all sums converge absolutely and uniformly on $|z| < \sqrt{R_1}$,
we conclude that $F_{R_2}(z)$ agrees with $F_{R_1}(z)$ there.
Letting $R \to \infty$, 
we obtain  an entire function $F_{\pi}(z)$ such that 
$$
F_{\pi}(z) = F_R(z) ~~~\mbox{when}~~~|z| \le \sqrt{R}.
$$
The maximum modulus bound \eqn{710} now establishes that $F_{\pi}(z)$
is an entire function of order one and exponential type at
most $\beta(\pi).$

Evaluating at integer points $n \in \ZZ$ and choosing $R=2n^2$
gives
\begin{eqnarray}~\label{713}
F_{\pi}(n) &= & 
 \sum_{|\rho| \le R}  \left(1 - (1 -\frac{1}{\rho})^n \right)
 -n \sigma_1(R, \pi) + 
\sum_{j=2}^{\infty} (-1)^j {\binom{n}{j}} \sigma_j(R, \pi)
\nonumber \\
&=&  \sum_{|\rho| \le R} \left(1 - (1 -\frac{1}{\rho})^n \right)
 -n \sigma_1(R, \pi) + 
\sum_{j=2}^{n} (-1)^j {\binom{n}{j}} \sigma_j(R, \pi).
\end{eqnarray}
The last term is a finite sum, so  we can expand all the 
terms in the first sum
$|\rho| < R$ and rearrange to obtain
\beql{714}
F_{\pi}(n) = -n \sigma_1( \pi) + \sum_{j=2}^{n} (-1)^j {\binom{n}{j}}
\sigma_j(\pi) = \lambda_{n}(\pi)
\eeq
using Lemma~\ref{Nle22} for the last equality. This proves property (i).

The symmetry property of zeros under $\rho \mapsto 1 - \bar{\rho}$
implies that $F_R(z)$ is real on the imaginary axis, for $|z| < R$.
It follows that $F_{\pi}(z)$ is real on the imaginary
axis, and  reflection principle then
gives the symmetry
\beql{715}
F_{\pi}(-x+ iy) = \overline{F_{\pi}(x+iy)}~~~\mbox{when}~~|z| < R.
\eeq
which is property (ii). 

We have constructed a specific function $F_{\pi}(z)$ 
which has properties (i), (ii). If  $\beta(\pi) \ge \pi$, then 
properties (i), (ii) do not determine the function uniquely
because we can add a real multiple of 
of $z\sin \pi z$ to the function while preserving  properties
(i) and (ii). 

(2) Now  suppose that the Riemann hypothesis holds for $\xi(s, \pi)$. 
Then we can write
\beql{716}
\rho = \frac{1}{2} + i \gamma = |\rho| e^{i \varphi_{\rho}}
\eeq
with $-\pi < \phi_{\rho} \le \pi$, and
\beql{716a}
\varphi_{\rho} := \tan^{-1} (2 \gamma),
\eeq
so that  $\phi_{\rho} \ge 0$
when $\gamma \ge 0$, and  $\varphi_{\rho}  < 0$ otherwise.
We have 
\beql{717}
1 - \frac{1}{\rho} = -\frac{\frac{1}{2} - i\gamma}{\frac{1}{2} + i\gamma}
= e^{i(-2\varphi_{\rho} \pm \pi)}
\eeq
with the sign $\pm$ chosen so that the angle falls in $(-\pi, \pi].$
Let us now set
\beql{717a}
\varphi_{\ast}(\rho) := \{ \begin{array}{cl}
-2\varphi_{\rho} + \pi & \mbox{if}~ \varphi_{\rho} > 0 
~\leftrightarrow \gamma >0\\
-2 \varphi_{\rho} - \pi & \mbox{if}~ \varphi_{\rho}< 0 
~\leftrightarrow~ \gamma < 0.
\end{array}
\eeq
The case $\varphi_{\rho}=0$, corresponding to $s = \frac{1}{2}$, is
problematic, since the sign $\pm \pi$ must have  a discontinuity there.
For the moment  we use the convention that the $+$ sign is chosen.
Consequently 
\beql{718}
|\log \left(1 - \frac{1}{\rho}\right)| = |i \varphi^{\ast}(\rho)| \le \pi
\eeq
and equality can occur only if $\varphi^{\ast}(\rho)= \pi$,
which requires that $\varphi_{\rho}=0$ and $\rho = \frac{1}{2}.$
It follows that if the Riemann hypothesis holds for $Z(\pi)$
then $\beta(\pi) \le  \pi$, 
and if there is no zero at $s= \frac{1}{2}$ then 
$ \beta(\pi) < \pi.$

For uniqueness, suppose first that $Z(\pi)$
includes  no zero at $s= \frac{1}{2}$.
It is well known 
(from the sampling theorem) that 
an entire function of order one and exponential type 
less than $\pi$ is
completely specified by its values at integer points.
Indeed,  the difference of two such functions would be an  entire function
of order one and exponential type $\tau$ less than $\pi$ vanishing at all
integer points. Since a function of exponential type $\tau$ has at most
$(\frac{2\tau}{\pi} + o(1)) R$ zeros in a disk of radius $R$, as $R \to \infty$,
it must be identically zero. In this case  property (i), 
together with a growth bound on the function, characterizes
the interpolation  function uniquely.
In the cases  where $Z(\pi)$ has zeros at $s=\frac{1}{2}$, we
use a convention to handle their contribution and
define $F_{\pi}(z)$ uniquely. If it has exactly 
$k$ zeros there, we  define their
contribution to be
$$
f_{\frac{1}{2}}(z) := \frac{k}{2} \left( e^{-i\pi z} + e^{i \pi z}\right).
$$
This function  is real on the imaginary axis, and has 
$f_{\frac{1}{2}}(z) = f_{\frac{1}{2}}(-z)$,
and it  corresponds
to  assigning half of the zeros at $s= \frac{1}{2}$ the
 argument $-\pi i$ and the other half argument $\pi i$.
Subtracting this contribution leaves an
interpolation  function having exponential type strictly
less than $\pi$, which is then uniquely specified by condition (i).

It remains to bound the size of the interpolation function
$F_{\pi}(z)$ on the real axis. Under the
Riemann hypothesis, the  bounds on the
Li coefficients in \S5 and \S6 lead one to expect an
upper  bound of shape $O\left((|x|+2) \log (|x|+2)\right),$
and we now show this holds.
We have
\beql{719}
F_{\pi}(x) = \sum_{\rho} \left( 1 - (1 - \frac{1}{\rho})^x \right)
\eeq
and we bound each term separately in absolute value.
For a zero $\rho$  with $\gamma >0$ we have
$$
\varphi_{\rho} = \frac{\pi}{2} - \frac{1}{2 \gamma} + O (\frac{1}{\gamma^2}).
$$
It follows that for $\gamma > 5|x|$ we have
\begin{eqnarray}
 1 - (1 - \frac{1}{\rho})^x &= &1 - e^{ix(-2 \varphi_{\rho} + \pi)}
\nonumber \\
&= & 1 - \exp \left(2ix(\frac{1}{2\gamma}+O(\frac{1}{\gamma^2})\right)  
\nonumber \\
&=& i \frac{ix} {\gamma} + O(\frac{x}{\gamma^2}). 
\end{eqnarray}
We find a similar expression when $\gamma <0$.

We now divide the sum \eqn{719} into parts $|\rho|< |x| + 2$
and $|\rho| \ge |x|+2.$ The former sum contributes 
$O\left((|x|+2) \log (|x|+2)\right)$ using the bounds on the
number of zeros in Theorem~\ref{th21}(4).
On the remaining range we pair the zeros with increasing positive
$\gamma$ against those with increasing negative $\gamma$, and
find that their imaginary parts cancel out by
Theorem ~\ref{th21}(4) to an error
$O(\frac{x \log (|\gamma|+2)}{\gamma^2})$ on each zero.
Now summing over the $O(\log T)$ zeros in an interval of
length one at height $T$, for $|x|+2 \le T \le \infty$
gives a total contribution of
$$
O \left( x (\log |x|+2)^2  \frac{1}{|x|+2} \right) = 
O\left((\log |x|+2)^2\right).
$$
The bound on the sum \eqn{719} is then
$$
|F_{\pi}(x)| \le C (|x|+2) \log(|x|+2),~~~-\infty < x < \infty,
$$
as asserted.
$~~~\bsq$ \\

We next consider the Fourier transform of $F_{\pi}(x)$
on the real line, regarded as a distribution. 
We suppose the Riemann hypothesis holds for $\pi$, and
consider the unique function $F_{\pi}(z)$ given
in Theorem~\ref{th71}(2).
The growth
bound in \eqn{N702a} implies that it is well-defined  as
a tempered distribution. Viewing \eqn{719} term-by-term
we obtain, a representation for this 
Fourier transform formally as a sum of delta functions
$$
\hat{F}_{\pi} (\eta) = 
\sum_{\rho \in Z(\pi)} 
\left( \delta_{0}(\eta) - \delta_{\varphi^{\ast}(\rho)}(\eta)
\right),
$$
in which $\delta_{\alpha}(\eta)$ denotes a delta function
centered at $\eta= \alpha$. Note that the individual terms
in the sum must be grouped as indicated to define
a continuous linear
functional on test function in the Schwartz space $\sS(\RR)$.
One deduces that  the support of this tempered distribution 
is real and lies in $[-\pi, \pi]$, and is a discrete set,
determined by the delta functions, except at the non-isolated
limit point $\eta=0$.

\paragraph{Remarks.}
(1) One can carry out a similar procedure to interpolate
the archimedean contributions  $S_{\infty}(n, \pi)$,
for $n \ge 1$ with an 
interpolation function
$F_{\pi, \infty}(z)$ that is  an entire function of
exponential type. 
One  must modify the proof to subtract
off the contribution of the term linear in $z$
from the
other terms in the formula of Lemma~\ref{le42},
to get a convergent formula.  The resulting
interpolation function is, formally,
\beql{731}
F_{\pi, \infty}(z) = C_3(\pi) + C_4(\pi) z + \sum_{j=1}^N \sum_{n=1}^\infty
\left(1 - \frac{2} {2n +\kappa_j(\pi)-1} z   - 
(1 - \frac{2}{2n +\kappa_j(\pi)-1})^{z} \right)
\eeq
for a certain constants $C_3(\pi), C_4(\pi)$. Here we have
$$
| 1 - \frac{2}{2n +\kappa_j(\pi)-1}| < 1~~~\mbox{for ~ all} 
~~ n \ge 1.
$$
The interpolation
function $F_{\pi, \infty}(z)$ will be an entire function of
exponential type. 
It  has polynomial growth on  the positive real axis 
(bounded by  $O((|x|+2) \log (|x|+2)$ for positive real
$x$),  but increases exponentially on the negative real axis.
Consequently its  Fourier transform can be interpreted as a 
distribution, but not as a tempered distribution.

(2) One can also obtain an entire function  
$F_{\pi,f}(z)$ of exponential type
interpolating $S_{f}(n, \pi)$,
given  as 
$$
F_{\pi,f}(z) := F_{\pi}(z) - F_{\pi, \infty}(z) - \delta(\pi).
$$
Again, even assuming the
Riemann hypothesis, its Fourier transform must be interpreted as a
distribution and not a tempered distribution.


(3)  Assume the Riemann hypothesis holds for $L(s, \pi)$,
and suppose that $\pi = \pi^{\vee}$ is self-dual. 
Then \eqn{717} yields 
\beql{720}
(1 - (1 -\frac{1}{\rho})^n) + (1 - (1 -\frac{1}{1- \rho})^n)
= 2 - 2 \cos n\varphi^{\ast}(\rho).
\eeq
Note that $\varphi_{\rho} \to \pm \frac{\pi}{2}$ as $\gamma  \to \pm\infty$,
and the $\pm \pi$ term is chosen of opposite sign, whence 
$-2 \varphi_{\rho} \pm \pi \to 0$ as $|\rho| \to \infty$.

%
%
%
\section{Concluding Remarks}
\hsp
(1) Instead of the $\xi$-function treated in \S2,
one may alternatively consider generalized
Li coefficients associated to  the function $\xi^{+}(s, \pi)$
obtained by
removing any zeros at $s= \frac{1}{2}$ from $\xi(s, \pi)$.
Here 
$$ 
\xi^{+}(s, \pi) := (s- \frac{1}{2})^{-e(\frac{1}{2}, \pi)} \xi(s, \pi),
$$
and we let $\lambda_n^{+}(\pi)$ denote the associated
Li coefficients. 
Use of this function has three positive features. 
First, it allows us to unambiguously  obtain the
functional equation in the form 
$$
\xi^{+}(s, \pi) = \xi^{+}(1-s, \pi)
$$
with the sign convention that $\xi^{+}(\frac{1}{2}) > 0$. 
Second, assuming RH, 
the interpolating function $F_{\pi}^{+}(z)$ 
for the generalized Li coefficients 
associated to $\xi^{+}(s, \pi)$ in \S7
will be an entire function of order $1$ and type strictly
less than $\pi$. In consequence it is uniquely determined
by the values $\{F_{\pi}^{+}(n, \pi) := \lambda_n^{+}( \pi), ~~ n \in \ZZ\}$.
Third, in terms of the Weil scalar product associated to
a hypothetical Hilbert-Polya operator, with eigenvalues
$\lambda = s^2 - \frac{1}{4}$,
removal of the zeros at $s= \frac{1}{2}$ would correspond  to
taking the orthogonal complement of the eigenspace with
eigenvalue $\lambda=0$. It is hoped that zeros at
$s= \frac{1}{2}$ 
will have some  arithmetic-geometric meaning, as
in the Birch-Swinnerton Dyer conjecture, and 
there may
well  be an arithmetico-geometric 
way to directly characterize this
eigenspace. \\

(2) There are analogues of the Li coefficients 
for automorphic $L$-functions in  the function
field case. We first note that 
the $\xi$-function for the trivial representation over a
function field $K$ in one variable
over a finite field $\FF_{q}$ can be taken to be 
$$
\xi(s, \pi_{triv,K})
= q^{-s}(1-q^{-s})(1 - q^{1-s})  Z_K(s),
$$
in which $Z_K(s)$ is the (completed) function
field zeta function, and is a polynomial in
$w= q^{-s}$. For a rational function field $K= \FF_q(T)$
we have $\xi(s, \pi_{triv, K})=1$.
All other automorphic $L$-functions are
polynomials in the variable $w$. 
The transformation $s = - \frac{z}{1-z}$ is still used
to define the Li coefficients in the function field case.
The definition \eqn{301} for the ``Weil scalar product''
applies, and the positivity of the 
Li coefficients is interpretable
in terms of this scalar product applied to  the 
 Li test functions  $G_n(s)$. The
asymptotics of the Li coefficients in the function field
case are different from those in the number
field case, since there are no archimedean places.
The main term in their asymptotics is $Cn$ rather
than the term $Cn\log n$ occurring in the number field
case.

To interpret the Li coefficients in terms of a 
function field ``explicit formula,''  one must
use a version in terms of  the $s$-variable,
noting that  $L$-functions are  singly periodic with period
$\frac{2\pi i}{\log q}$. 
The standard  function field ``explicit formula'' is 
generally stated in terms of the 
variable $w$, and it 
has an algebraic geometry interpretation,
related to intersection theory.  The allowed
test functions in this formula
(related to divisors) are 
Laurent polynomials in $w= q^{-s}$ with integer
coefficients, cf. Haran~\cite{Ha91}. 
The  Li test functions 
$G_n(s)$ are not functions of $q^{-s}$, and must be
$q$-periodized using 
\beql{803}
P_{q}(G_n)(s) := \sum_{n \in \ZZ} G_n(s+ \frac{2\pi i n}{\log q}).
\eeq 
to be viewed as function field test functions. 
This sum \eqn{803} is conditionally convergent and is regularized
as a limit as $T \to \infty$ of $\sum_{-T < n < T}$.
It would be interesting to see if this framework is  compatible with a
framework of Deninger~\cite{De92}, \cite[Sect. 4]{De94a}.
Another feature of the function field
case is that some function field $L$-functions
have multiple zeros on the critical line at positions
other than $s \equiv  \frac{1}{2}~ (\bmod \frac{2\pi i}{\log~q})$.\\

(3) Under the Riemann hypothesis, the ``explicit formula'' decomposition of 
the Li coefficients into $S_{\infty}(n, \pi)$
and $S_f(n, \pi)$ 
reveals that 
the dominant contribution to the asymptotics of the Li coefficients
comes from  the archimedean terms, which 
correspond in \eqn{731}
to the ``trivial zeros'' of $L(s, \pi)$.
This contrasts with
the definition \eqn{101} of the Li coefficients \eqn{101},
which is  a sum over  the non-trivial zeros of $L(s, \pi)$,
and does not include the ``trivial zeros''.
In effect certain asymptotics of the nontrivial
zeros are  described in terms of the  ``trivial zeros.'' 


%
%
%
\newpage
\section{Appendix: ``Explicit Formula'' and Weil's Quadratic Functional}
\hsp

For simplicity  we treat here only the
case of the trivial representation $\pi_{triv}$ 
on $GL(1)$. In that case, there are two versions of the
 ``explicit formula'' of prime number theory,
the trace form and the covariance form, described
below.
For all other automorphic representations $\pi$ there
is essentially only one form of the ``explicit formula.''

Consider first a vector space $\sA_{\delta}$
 of test functions $f: \RR_{>0} \to \CC$, defined 
as follows.
Associated to a test function is its Mellin transform
$$
\hat{f}(s) := \int_0^{\infty} f(x) x^s \frac{dx}{x}.
$$
The vector space $\sA_{\delta}$ consists
of those functions $f(x)$ whose  Mellin
transform is analytic in the  strip 
$\frac{1}{2} - \delta < \Re(s) < \frac{1}{2} + \delta$,
and extends continuously to  the boundary of the strip, 
We put a metric on test functions given by
the uniform norm on the closed strip,
$$
d(f_1, f_2) := sup_{s \in S_{\delta}} |f(s) - g(s)|.
$$
With this topology $\sA_{\delta}$ is a complete metric
space.  
We define the involution 
\beql{901a}
\tilde{f}(x) := \frac{1}{x} f( \frac{1}{x}),
\eeq
whose effect on Mellin transforms is:
\beql{902a}
\hat{\tilde{f}}(s) = \hat{f}(1-s).
\eeq
We define the space $\hat{\sA_{\delta}}$ to be the set of Mellin
transforms of the elements of $\sA_{\delta}$, regarded as 
analytic functions in the specified strip. Note that the space $\sA$
in  \S3 is given by 
$$
\sA = \cup_{\delta< 1/2} \hat{\sA}_{\delta}.
$$

We define the {\em Weil distribution functional}
\beql{903a}
W[f] := \sum_{ \{\rho: \xi(\rho) = 0 \}}{}^{'} \hat{f}(\rho),
\eeq
in which ${}^{'}$ means that the 
(possibly conditionally convergent) sum is interpreted
as $\lim_{T \to \infty} \sum_{ |\rho| \le T} $.
Under suitable conditions $W[f]$ is a continuous linear
functional on the allowed set of test functions $\sA_{\delta}$,
for example when $\delta > \frac{1}{2}.$ ( If the Riemann hypothesis
is assumed one can take any $\delta > 0.$)
We then   define the {\em trace functional}
\beql{904a}
T[f] := \hat{f}(0)  - W[f] +\hat{f}(1)
\eeq
The quantity $T[f]$ is sometimes called the ``spectral side''
of the ``explicit formula'' of prime number theory.

The  ``explicit formula'' in {\em trace form} is an
formula  for $T[f]$ 
for a suitable set of test functions $f$,
taking the shape 
\beql{trace}
T[f] = \sum_{\nu} W_{\nu}(f),
\eeq
in which $W_{\nu}(f)$ is a contribution associated to each
(non-archimedean or archimedean) place $\nu$ of the given
field $K$.
The right hand side of \eqn{trace} is sometimes called the ``arithmetic-
geometric'' side of the explicit formula. We shall not be
concerned with the exact form of the 
arithmetic-geometric  side here, the 
individual terms $W_{\nu}(f)$  of which  can be expressed
in various interesting ways, see Burnol \cite{Bu98}, Haran \cite{Ha90}.
For general automorphic representations
$\pi$  the terms $W_{\nu}(f)$ represent contributions from 
individual terms in the Euler product factorization of
$\xi(s, \pi)$.
A version of the "explicit formula" in the trace
formulation is given in Patterson  ~\cite[Sect. 3.6]{Pa88}.
In order to get an unconditional result one must use
a test function space 
$\sA_{\delta}$ with $\delta > \frac{1}{2}.$
A version of the "explicit formula"  for principal
automorphic $L$-functions over $GL(n)$ appears
in Rudick and Sarnak~\cite[Prop. 2.1]{RS96}.

One hope is  that the
``explicit formula''  \eqn{trace}
might be interpretable as a Lefschetz trace formula
coming from a dynamical system acting on an (unknown)
``geometric'' object,  in which
$T[f]$ is to be viewed as a distributional trace with
test functions $f$ living on cohomology groups of the geometric
object, and the arithmetic side giving
data arising from fixed points or periodic orbits on the
geometric object. This viewpoint is taken in  
Deninger~\cite{De93}, \cite{De94} \cite{De98}, \cite{De02}.

Weil originally formulated  the explicit formula 
in the {\em covariance form} given by 
\beql{covariance}
W[f] = -\sum_{\nu} W_{\nu}(f) + W_{0}(f) + W_{1}(f).
\eeq
This is a rearrangement of the terms in the ``trace form'' equality,
in which we define 
$$
W_{0}(f) := \hat{f}(0) = \int_0^{\infty} \frac{1}{x}f(x) \frac{dx}{x}, 
$$
and
$$
W_{1}(f) := \hat{f}(1) =\int_0^{\infty} f(x) \frac{dx}{x}.
$$
Attached to this quantity Weil introduced the {\em scalar product}
\beql{905a}
\langle f, g \rangle_{\sW} := W[ f * \tilde{\bar{g}}]
= \sum_{\rho} \hat{f}(\rho) \overline{\hat{g}(1-\bar{\rho})}
=  \sum_{\rho} \hat{f}(\rho) \hat{\hat{\bar{g}}}(1 - \rho)   .
\eeq
Weil's criterion for the Riemann hypothesis is that this
scalar product be positive semidefinite on a suitable
space of test functions, which could be taken to be
$\sA_{\delta}$ for any $\delta > \frac{1}{2}$, for example.
This scalar product will be indefinite  if the
Riemann hypothesis does not hold.
We may transport this scalar product  forward to
the Mellin-transformed space $\hat{\sA}$ as
$$
\langle \hat{f}, \hat{g} \rangle_{\sW} := \langle f, g \rangle_{\sW}.
$$
The Weil functional $W[f]$ has a natural  interpretation in the function
field case, see Haran~\cite{Ha91}.

One feature of the covariance form \eqn{covariance}
of the ``explicit formula'' compared to the 
trace form is that one can make
sense of the functional $W[f]$  on  a larger set of test functions 
than those allowed in the trace form, permitting 
test functions that have  singularities at $s=0$ and $s=1.$
However the  right hand side of \eqn{covariance}
must then  be redefined
as a limit as $T \to \infty$, with each local term 
computed using cutoffs at $x= \frac{1}{T}$ near zero and at $x= T$
near $\infty$.
This was done in the computations in \cite{BL99}.
This is a "regularization" of the right side of \eqn{covariance}
because the
two terms $W_{0}(f)$ and $W_1(f)$ diverge as the
cutoff parameter $T \to \infty$. Note that these two terms
are present only for the trivial representation $\pi_{triv}$
on GL(1), and
are absent from the analogous formula for all other 
cuspidal automorphic representations $\pi$.
In these other cases the analoguous formulas have  $T[f] = - W[f]$,
so the ``trace form'' and ``covariance form'' essentially coincide.

The vector space $\sL$ of  Li test functions makes
sense for this extended covariance form
of the "explicit formula" with a cutoff parameter, 
and not for the trace form. 
As noted in \S3, it
consists exclusively of  
rational functions which
have poles either at $s=0$ or $s=1$, or both. In
consequence the trace function  $T[f]$ is undefined for every 
Li test function. 
For example, if 
$g_n(x)$ denotes the test function on $\RR_{>0}$
whose Mellin transform 
$$
G_n(s) \hat{g}_n(s) = \int_{0}^\infty g_n(x) x^{s-1} dx,
$$
gave the Li test functions 
$G_n(s)$
in \eqn{303}
then 
$\hat{g_n}(0)$ is infinite, for each $n \ge 1$. 
In \cite[Lemma 2]{BL99}
the functions $g_n(x)$
were explicitly determined, as 
$$
g_n(x) = \{ \begin{array}{cl}
P_n(\log x) & \mbox{if}~~ 0 < x < 1 \\
\frac{n}{2} & \mbox{if}~~ x = 1 \\
0 &  \mbox{if}~~ x > 1.
\end{array}
$$
in which 
$$
P_n(x) :=\sum_{j=1}^n {\binom{n}{j}} \frac{x^{j-1}}{(j-1)!}.
$$
We note  that $P_{n}(x)= L_{n}^{{\alpha}}(x)$ is a Laguerre polynomial
with $\alpha=1$.
The Weil scalar
product of $g_n$ with other functions  remains well-defined
in the covariance form of the ``explicit formula'' when
a cutoff version is used. 
In terms of the cutoff version of the definition,
the unbounded  contribution from $\hat{g_n}(0)$
above as $T \to \infty$ is offset by a corresponding
divergence coming from the finite primes in
the ``explicit formula.''
The Weil distribution functional  $W[f]$ is well-defined
for all $f \in \sL$, as is the Weil scalar product, whether or
not the Riemann hypothesis holds.

\end{document}